\documentclass[11pt]{article}

\usepackage[numbers,sort&compress]{natbib} %\cite{1,5-8,11}
\usepackage{color}
\pdfoutput=1

\usepackage{amsfonts}
\usepackage{mathrsfs} %book, letter, report, % none for 10pt
\usepackage{amsfonts,amssymb}
\usepackage{exscale}

\usepackage{comment}

\usepackage{amsgen,amsmath,amstext,amsbsy,amsopn,amsfonts,amssymb,pifont, epsfig}
\usepackage{graphicx}

\newtheorem{thm}{Theorem}[section]
\newtheorem{lemma}[thm]{Lemma}

\newtheorem{remark}[thm]{Remark}

\newcommand{\be}{\begin{equation}}
\newcommand{\ee}{\end{equation}}
\newcommand{\bea}{\begin{eqnarray}}
\newcommand{\eea}{\end{eqnarray}}
\newcommand{\beaa}{\begin{eqnarray*}}
\newcommand{\eeaa}{\end{eqnarray*}}
\newcommand{\bei}{\begin{itemize}}
\newcommand{\eei}{\end{itemize}}
\newcommand{\bee}{\begin{enumerate}}
\newcommand{\eee}{\end{enumerate}}
\newcommand{\bi}{\begin{itemize}}
\newcommand{\ei}{\end{itemize}}
\newcommand{\beq}{\begin{eqnarray*}}
\newcommand{\eeq}{\end{eqnarray*}}
\newcommand{\beqn}{\begin{eqnarray}}
\newcommand{\eeqn}{\end{eqnarray}}

\newcommand{\pf}{\noindent {\bf Proof \ }}

\newcommand{\ignore}[1]{}{}

\def\FF{\mathcal{F}}

\def\OO{\mathcal O}

\def\qed{{\hfill $\Box$ \bigskip}}

\def\P{{\mathbb P}}

\def\R{{\mathbb R}}

\newcommand{\vp}{\varepsilon}
\def\<{\langle}
\def\>{\rangle}
\def\wh{\widehat}
\def\wt{\widetilde}

\def\boo{\bar{\mathcal O}}
\def\po{\partial{\mathcal O}}
\def\R{\mathbb R}
\def\T{\mathbb T}

\newtheorem{assumption}{Assumption}

\numberwithin{equation}{section}

\usepackage{color}
\definecolor{webgreen}{rgb}{0,.5,0}
\definecolor{webbrown}{rgb}{.8,0,0}
\definecolor{emphcolor}{rgb}{0.5,0.95,0.95}
\usepackage{hyperref}
\hypersetup{%
          colorlinks=true,
          linkcolor=webbrown,
          filecolor=webbrown,
          citecolor=webgreen,
          urlcolor=blue,
          breaklinks=true}

\begin{document}

%\includegraphics[width=18.5cm, higth=10cm]{file.jpg}
%\title{\bf The third boundary value problems of semilinear elliptic PDEs and BSDEs with singular coefficients }
\title{\bf Homogenization of semilinear parabolic PDEs with the third boundary conditions}

\author{Junxia Duan, Jun Peng$\thanks{The corresponding author. Email: pengjun2015@csu.edu.cn, duanjunxia0802@163.com}$\\ {\small School of Mathematics and Statistics, Central South University, Changsha, Hunan 410083, China}}

\maketitle

\begin{abstract}
  In this paper, we study the homogenization of the third boundary value problem for semilinear parabolic PDEs with rapidly oscillating periodic coefficients in the weak sense. Our method is entirely probabilistic, and builds upon the work of \cite{Tanaka} and \cite{B Hu}. Backward stochastic differential equations with singular coefficients play an important role in our approach. %semilinear PDEs with the third boundary condition, periodic coefficients and highly oscillating drift and
\end{abstract}

\smallskip\noindent
{\bf Keywords}: Homogenization; Weak solution; Third boundary value problem; Backward stochastic differential equations

\smallskip\noindent
{\bf Mathematics Subject Classification (2020)}: 60H30, 35B27, 35K40
%35  PDEs
%35B27 Homogenization in context of PDEs; PDEs in media with periodic structure;    35D30 Weak solutions to PDEs;     35G60 Boundary value problems for systems of nonlinear higher-order PDEs;       35K40 Second-order parabolic systems;        35K61 Nonlinear initial, boundary and initial-boundary value problems for nonlinear parabolic equations

%60 Probability theory and stochastic processes
%60H30  	Applications of stochastic analysis (to PDEs, etc.)  ;  	60J46  	Dirichlet form methods in Markov processes;     60J50  	Boundary theory for Markov processes;     	60J55  	Local time and additive functionals;    	   60J60  	Diffusion processes [See also 58J65]

%58J32  	Boundary value problems on manifolds

%31C25, Dirichlet form,
%60H15,  SPDE(aspects of stochastic analysis,
%60J60 Diffusion processes,
%Schrodinger and Feynman-Kac semigroup

\bigskip

\section{ Introduction }

%partial differential equations (PDEs) \quad
%study the limit behaviour, as $ \varepsilon \to 0$, of the solution $ u^\varepsilon$
In the study of porous media, composite materials and other systems of physics and mechanics, they frequently involve the boundary value problems with periodic structures. %The aim of homogenization theory is to establish the macroscopic rigorous characterizations of the microscopically heterogeneous media.
The process of establishing the macroscopic rigorous characterizations of such microscopically systems is called homogenization. It has been a highly active research area in mathematics for a long time, and a vast literature exists on this topic. See \cite{B Hu,Castell,G P,Lejay a,Lejay b,Pardoux} and references therein.

%The commonly used methods can be classified into analytical methods and probabilistic ones.
Homogenization theory has motivated the development of various notions of weak convergence in analysis. Such convergence can be better understood from the direction of probabilistic interpretation of the equation.
Generally speaking, the probabilistic method begins with representing the quantities of interest by means of stochastic processes, and then attempts to prove convergence in laws of these processes. Hence it is also known as the averaging principle.

The goal of this paper %In this paper, our primary interest
is to use a probabilistic approach to study the homogenization %limit behaviour, as $ \varepsilon \to 0$, of the solution $ u^\varepsilon$
of the following third boundary value problem for semilinear parabolic PDEs with rapidly oscillating periodic coefficients
%of the following form
\begin{equation}\label{e:1.1}
\begin{cases}
\frac{ \partial u^\varepsilon}{ \partial t} (t,x)+ L^\varepsilon u^\varepsilon(t,x) +  f(x,u^\varepsilon(t,x),\nabla u^\varepsilon(t,x) )=0   & (t,x) \in [0,T) \times \OO,  %x\in\OO, ~t\in [0,T),
\cr
%\gamma(\frac{x}{\varepsilon} )  \cdot \nabla u^\varepsilon(t,x)
\frac12 \frac{ \partial u^\varepsilon }{\partial \upsilon^\vp } (t,x)
 +c(x/ \vp ) u^\varepsilon(t,x)=0   & (t,x) \in [0,T) \times \po, %x\in \po, ~t\in [0,T),
\cr
u^\varepsilon(T,x)= g(x) & x\in \boo,
\end{cases}
\end{equation}
where $\OO$ is a bounded convex domain %with {\color{red}%$C^3$
%smooth} boundary
in $\R^d  $ ($d\geq 2$) and the boundary $ \po$ is said to be of class $ C_b^1$, that is there exists a function $\Psi\in C_b^1 ( \R^d, \R)$ such that
\begin{equation}\label{e:o}
\OO= \{ x\in \R^d: \Psi (x) >0 \}, \quad \inf_{x\in \po} | \nabla \Psi (x)| >0.
\end{equation}
Here
$
\frac{ \partial u^\varepsilon }{\partial \upsilon^\vp }
=  \vec{n}_i(x)  a_{ij} \big( \frac{x}{ \vp}  \big)
\frac{\partial u^\vp}{ \partial x_j}
$
denotes the conormal derivative associated with $ L^\vp$,
%$\gamma(x/\varepsilon  ) = \frac12  A(x/\varepsilon )  \vec{n}(x)$, where
and $\vec{n}$ is the inward unit normal to $\po$.
The family divergence-form operators $L^\varepsilon$ are given by
\begin{eqnarray}
L^\varepsilon &:=& \frac12 \nabla \cdot \left(A\big( \cdot/\varepsilon \big) \,  \nabla \right) + \frac{1}{\varepsilon} b\big( \cdot/\varepsilon \big) \cdot \nabla
\nonumber \\
\label{e:1.2}
&=&\frac12 \sum_{i,j=1}^d
\frac{\partial}{ \partial x_i}
\left(a_{ij} ( \cdot / \varepsilon)  \frac{\partial}{ \partial x_j} \right)  + \frac{1}{\varepsilon}  \sum_{i=1}^d
b_i (\cdot /\varepsilon)\frac{\partial}{ \partial x_i},
\end{eqnarray}
where the coefficients $ a$ and $b$ are periodic.
%For notational simplicity, we use Einstein¡¯s convention that the repeated indices in a product will be summed automatically.

%\begin{assumption}\label{h:1}
$A(x)= ( a_{ij}(x) )_{1\leq i,j \leq d} : \R^d \to \R^d  \otimes  \R^d $ is a smooth%measurable
, symmetric matrix-valued function satisfying the uniformly elliptic condition
\begin{equation}\label{e:1.3}
\lambda^{-1} I_{d\times d} \leq A(\cdot) \leq \lambda I_{d\times d},
\end{equation}
for some constant $\lambda \geq 1$.
%In the sequel, Moreover,
%\end{assumption}
The function $c:  \R^d \to [-\alpha,0]$ is non-positive for some constant $ \alpha \in [0,\infty)$ and periodic with respect to the orthonormal basis $ \{(e_1, \dots, e_d ) \}$ of $ \R^d$. i.e. %More precisely, $ c$ satisfies
\begin{eqnarray*}
 c\in C_b(\R^d)   \qquad \quad\hbox{and}~
\qquad \quad c( x+ e_i )= c(x)\leq 0,
\end{eqnarray*}
for $i= 1,\dots, d $.
%non-positive bounded and defined on the boundary $\po$.
%We can see Assumption \ref{h:3} in Section \ref{S:2} for the details of the definitions and conditions on $ f$.
Assumptions \ref{h:2}-\ref{h:4} are made for the coefficients and will be listed in Section \ref{S:2}. %the sequel
The main result of this paper is the following theorem, which we will prove at the end of Section \ref{S:5}.

\begin{thm}\label{T:1.1}
Let Assumptions \ref{h:2}-\ref{h:4} hold, then
the semilinear PDE \eqref{e:1.1} has a unique weak solution $ u^\varepsilon$ for each $ \vp >0$. Moreover under \eqref{e:3.5} and \eqref{e:n}, we have for each $ t\in [0,T], x\in \boo$,
%\begin{equation}\label{e:1.4}
$$
u^\vp (t,x) \to u^0(t,x), \qquad \qquad  ~\hbox{as}~\vp \to 0,
%\end{equation}
$$
where $ u^0$ satisfies the limit semilinear PDE
\begin{equation}\label{e:1.5}
\begin{cases}
\frac{ \partial u^0}{ \partial t} (t,x)+ \overline{L} u^0(t,x) +  \bar{f}(x,u^0(t,x),\nabla u^0(t,x) )=0   & (t,x) \in [0,T) \times \OO,  %x\in\OO, ~t\in [0,T),
\cr
\frac12 \frac{\partial u^0}{ \partial  \upsilon^0    }  (t,x)
%\partial_{ \bar{\Gamma}} u^0(t,x)
+\bar{C} u^0(t,x)=0   & (t,x) \in [0,T) \times \po,  %x\in \po, ~t\in [0,T),
\cr
u^0(T,x)= g(x) & x\in\boo.
\end{cases}
\end{equation}
Here%In \eqref{e:1.5},
\begin{equation}\label{e:1.6}
\overline{L} := \frac12  \sum_{i,j=1}^d  \frac{\partial }{ \partial x_i}    \left( \bar{a}_{ij}  \frac{ \partial}{ \partial x_j}    \right),
\end{equation}
$\bar{f} $ is given by \eqref{e:5.1} and
$\frac{\partial u^0}{ \partial   \upsilon^0    }  $ is the conormal derivative associated with $ \overline{L} $.
%and $ \upsilon^0$ coincides with the constant vector $  \bar{\Gamma}=  (\bar{\Gamma}_i)_{1\leq i\leq d} $ given by ......
Moreover, the homogenized (or effective) coefficients $\bar{A}%\bar{a}_{ij}
=( \bar{a}_{ij}  )_{1\leq i,j \leq d}$ and $ \bar{C}$
are constant, which are given by \eqref{e:a} and \eqref{e:5.2} respectively.
\end{thm}

For the linear case (i.e. $f=0$ in \eqref{e:1.1}),
the original probabilistic approach to the homogenization of %local linear
the second order parabolic partial differential operators
is presented in Chapter 3 of \cite{Bensoussan}, which is based on the thoughts of Freidlin \cite{Freidlin}, that is finding harmonic functions (the so-called ``auxiliary" problem). Hence the problem can reduce to transforming the underlying Markov process into a martingale, whose quadratic variation has a deterministic limit by the ergodic theorem.
%ergodic theorem, the Feynman-Kac formula and the functional central limit theorem.
By now, this method has been extended in various directions, such as in the case of periodic coefficients \cite{Lejay a,Lejay b,Pardoux,Tanaka} and the case when the coefficients are stationary random fields \cite{Castell,Lejay b,Papanicolaou}.
%a lot of homogenization problems have been investigated for various elliptic and parabolic differential equations.
But for the nonlinear case (i.e. $f \neq 0$ in \eqref{e:1.1}), such PDEs are generally very hard to solve.
Since the notion of nonlinear backward stochastic differential equations (BSDEs) was first introduced by Pardoux and Peng~\cite{PP}, BSDEs have been used effectively to solve certain PDEs with nonlinear terms (see \cite{BPS,Ouknine Pa,Rozkosz} and references therein).
In view of the connection between them, the probabilistic tool has also been widely used in order to prove homogenization results for certain classes of nonlinear PDEs. %with highly oscillating drift.
When rapidly oscillating nonlinear terms are of the type $ f( t,x/\vp, u^\vp)+ \nabla u^\vp \, \wh f( t,x/\vp, u^\vp) $ (see \cite{Pardoux,Lejay b}) %is not depending of the gradient (more precisely, ),
or the type $ f(x, u^\vp) + \wh f(x, u^\vp) \, \| \nabla u^\vp\|^2$ (see \cite{G P}), it successfully proves convergence in laws of the stochastic processes involved with the help of BSDEs. However, for %the case of
more general nonlinearity in the gradient, we need another method in \cite{B Hu}, which is to exploit the stability results of BSDEs. This strategy has also been employed in the homogenization of random PDEs \cite{Castell}.

%Without the boundary conditions,

In this paper, we stress that the solution we considered for PDEs are weak solutions, not viscosity solutions.
In addition, we are interested in homogenization problems for which it is necessary to identify both the homogenized PDEs and the homogenized boundary conditions.
There are few works dealing with homogenization of PDEs with boundary conditions by probabilistic methods (see \cite{Bensoussan,Lejay b,PS} for Dirichlet boundary conditions and \cite{Chen Wu,Ouknine Pa} for Neumann boundary conditions %in half-space
).
We mainly consider PDEs with the third boundary boundary conditions, which will give rise to several difficulties both in analytic and probabilistic aspects. %We will describe them in detail below.
Firstly, due to the asymmetry of operator $L$ in the auxiliary problems \eqref{e:3.1}, we cannot expect an analogue of Meyers'result \cite[Theorem 2]{Meyer} to hold for every $p>2$.
%\cite[Theorem 1]{Groger} shows, however, that under weak hypotheses \eqref{e:3.5} we are able to prove it holds for some $p>2$
We are able to prove, however, that under weak hypotheses \eqref{e:3.5} it continues to hold for some $p>2$ by \cite[Theorem 1]{Groger}. Hence the $ L^p $-estimates $( p>2)$ for gradients of solutions can be established in Theorem \ref{T:3.2}, which plays an important role in the proof of homogenization in Subsection \ref{S:5.3}.
Secondly, for the property \eqref{e:5.6} of process $ X^\vp$, the method of \cite[Lemma 3.9.2]{Bensoussan} is no longer applicable because of the appearance of the local time term and the highly oscillating term $ \vp^{-1}  \wt b ( X^\vp /\vp) \, ds$ in \eqref{e:2.1}. Here we are highly motivated by \cite[Lemma 6.1]{Tanaka} and flexibly use the property of the local time.%\cite[Lemma 4.1]{B Hu}.
%The terms $ cu^\vp $ at the boundary conditions give rise to several difficulties both in analytic and probabilistic aspects. %This will give rise to several difficulties both in analytic and probabilistic aspects.
%. Inspired by \cite[Theorem 1]{Groger} and \cite[Theorem 2]{Meyer},

The organization of the paper is as follows.
In Section \ref{S:2}, we define a $ \T^d$-valued Markov process with generator $ L$ (defined by \eqref{e:L}) having periodic coefficients and a $ \T^{d-1}$-valued Markov process on the boundary with generator $ H_{\gamma}$ (defined by \eqref{e:2.3}). Then consider
the invariant measure of these two processes respectively. %: $  m(x) dx$ and $ \wt m(dx) $. %
We finally state the assumptions made throughout the paper.
In Section \ref{S:3}, we deal with the periodic solutions of the auxiliary problems \eqref{e:3.1} and obtain
the $ L^p $-estimates (some $ p>2$) of the gradient of the solution (Theorem \ref{T:3.2}). %are established in Theorem \ref{T:3.2}. %using the results of Gr\"{o}ger \cite{Groger}
In Section \ref{S:4}, we mainly study the existence and uniqueness of a class of BSDEs with singular coefficients (Theorem \ref{T:4.1}). %prove the existence and uniqueness of an $L^2$-solution $ (Y,Z)$ of the following BSDEs with
In Section \ref{S:5}, homogenization result for the problem \eqref{e:1.1} is proved.
%the third  boundary value problem of parabolic PDEs

We use the following notation in this paper.
For a matrix $\sigma$, its transpose and Hilbert-Schmidt norm are expressed by $\sigma^*$ and $\| \sigma \|= (\sum_{ij} \sigma_{ij}^2)^{\frac12}$. Denote by $|x|$ the Euclidean norm of $x$ in $\R^d$ and by $\< x,y\>$ the inner product of $x,y\in \R^d$.
%Define
%$ f_\vp (x):= f(x/\vp )$ (or $ F_\vp (x,y):= F(x/\vp, y  )$)
%for a function $ f$ on $ \R^d$ (or $ F$ on $ \R^d \times \R^d $ ).    %then the boundary condition in \eqref{e:1.1} can be written as
%$$
%\frac{ \partial u^\varepsilon }{\partial \upsilon^\vp } (t,x)
%+ c_\vp(x) u^\vp (t,x)=0, \qquad x\in \po, ~t\in [0,T).
%$$
The torus $\T^d:= \R^d/ {\mathbb Z}^d $ will be used frequently and we shall always identify the periodic function on $ \R^d$ of period $1$ with its restriction on the torus $\T^d$.
Denote by $ L^2( \T^d)$ and $ H^1( \T^d)$ the space of functions locally in $ L^2( \R^d)$ and $ H^1( \R^d)$ which are $ \T^d$-periodic.
%From the compactness of the torus $ \T^d$, it yields that the injection from the space $ H^1( \T^d)$ to $ L^2( \T^d)$ is compact.
Thanks to the compactness of the torus $ \T^d$, the injection from the space $ H^1( \T^d)$ to $ L^2( \T^d)$ is compact.

\section{General assumptions and preliminaries }\label{S:2}
This section is devoted to finding a $ \T^d$-valued Markov process with generator $ L$ (defined by \eqref{e:L}) having periodic coefficients and a $ \T^{d-1}$-valued Markov process on the boundary with generator $ H_{\gamma}$ (defined by \eqref{e:2.3}).
Then consider the invariant measure of these two processes respectively, which will be used to define the homogenized coefficients %in the proof of homogenization
in Section \ref{S:5}.
In the end, we state the assumptions more precisely made on the systems \eqref{e:1.1}.

%Let $(  X_x(t), \P_x, %\theta_t,
% x\in \OO )$ be the
Given $ \vp >0, x\in \OO$, the differential operator $ L^\vp$ inside $ \OO$ together with the Neumann boundary condition
$ \< \gamma (x / {\vp} )  , \nabla \cdot \> := \<A( x/ \vp) \vec{n}   , \nabla  \cdot\>=0 $
on $ \po$ determines a unique reflecting diffusion process $(  X^{\vp} (s), \P_{t,x}^{\vp}, %\theta_t,
t\in [0,T],x\in \OO )$ starting from $ x$ at time $ t$.
$ \mathcal{F}^{X^\vp}$ is the minimal admissible filtration generated by $ X^{\vp}$.
Set $\wt b:=(\wt b_1, \cdots, \wt b_d)^T $, where $\wt b_i= \frac12 \sum_{j=1}^d  \frac{\partial a_{ij}}{\partial x_j}   + b_i$.
Then by \cite{Lions},
%In general, $X^{\vp} $ is not a semimartingale, but
it has the following decomposition
\begin{eqnarray}\label{e:2.1}
X_s^{\vp} = x+ M_s^\vp +  \frac{1}{\vp} \int_t^s   \wt b   \left( X_r^{\vp}  /\vp \right)    \, dr + \int_t^s \gamma \left( X_r^{\vp}  /\vp  \right) \, d K_r^\vp, \quad 0\leq t\leq s\leq T, \qquad
\P_{t,x}^{\vp}-a.s.,
\end{eqnarray}
where $ M^\vp$ is a martingale additive functional of $X^{\vp} $ with quadratic cross-variation
$$ d \ll M^{\vp,i}, M^{\vp,j} \gg_s = a_{ij} ( X_s^{\vp} / \vp)  \, ds ,$$
and $K_s^\vp =\int_t^s I_{ \{   X_r^{\vp} \in \po   \} }  \,dK_r^\vp$ is the boundary local time of $X^{\vp} $.

%Next
Via the canonical quotient map $ \pi: \R^d \to \R^d / {\mathbb Z}^d$, we can define a $ \T^d$-valued Markov process with generator $ L$ (defined by \eqref{e:L}) having periodic coefficients. Meanwhile, we also want to find a $ \T^{d-1}$-valued Markov process on the boundary with generator $ H_{\gamma}$ (defined by \eqref{e:2.3}) by using the local time serves as a time change function.
%Denote by $m$ and $ \tilde{m}$
Now consider
the invariant measure of these two processes respectively.
\medskip

{\bf The invariant measure $ \bf m(x) dx$ :}
Let %$ L$ be the operator
\begin{equation}\label{e:L}
L:= \frac12  \sum_{i,j=1}^d   \frac{\partial}{ \partial x_i}
\left(a_{ij}   \frac{\partial}{ \partial x_j} \right)  + \sum_{i=1}^d
b_i \frac{\partial}{ \partial x_i},  %, \qquad
%D(L)= \big\{ f\in H^1(\T^d)  : Lf \in L^2 (\T^d)  \big\}.
\end{equation}
%According to \cite{Bensoussan, Lejay a}, we can define a Markov process $ \overline{X}$ corresponding to the projection on the torus $ \T^d$ of the process generated by the divergence-form operator $ L$ having periodic coefficients.
then it is well-known that the divergence-form operator $ L$ generates a Markov process on $ \R^d $.
By mapping all trajectories of this processes on $ \R^d $ to the torus $ \T^d$, %via the canonical quotient map $ \pi: \R^d \to \R^d / {\mathbb Z}^d$,
we can define a Markov process $ \overline{X}$, which is $ \T^d$-valued and generated by the operator $ L$ having periodic coefficients
(see \cite[Section 3.3.2]{Bensoussan} or \cite[Section 3]{Lejay a} for details).
In view of the compactness of the torus $ \T^d$, the process $ \overline{X}$ is ergodic.

Moreover,
from the maximum principle and the $ \T^d$-periodicity of functions in $ L^2 ( \T^d)$, it follows that any solution to $ L u=0$ is constant. Hence by Fredholm alternative theorem, there exists a unique solution to $ L^* m =0$ such that $ \int_{ \T^d} m(x) \, dx =1$.
According to \cite{Bensoussan, JKO},
we can obtain that $ m$ is positive, continuous and in fact the density of the invariant measure of the process $ \overline{X}$.
%the unique invariant probability of the process $ \overline{X}$.
\medskip

{\bf The invariant measure $ \bf \widetilde{m}(dx)$ :}
Given a function $ \varphi \in C(\po)$ with bounded partial derivatives of order $\leq 2$
and consider the problem%denote by $ \tilde{u}$ the solution of the problem
\begin{equation}\label{e:2.2}
\begin{cases}
L\tilde{u}(x)=0   &\hbox{in } \OO,   \cr
\tilde{u}(x)= \varphi(x) & \hbox{on }\po,
\end{cases}
\end{equation}
where $ L= \frac12 \nabla \cdot( A\nabla) +b \cdot \nabla  $.
According to \cite[Theorem 1.1]{Chen Zhao},
$
\tilde{u} (x)= E_x [ \varphi ( X_{\tau (\OO)} )  ]
$
is the unique continuous weak solution to Eq.\eqref{e:2.2},
where $ X$ is the continuous diffusion process generated by the operator $ L$ and $\tau (\OO) $ is the first exit time from $ \OO$, that is $\tau (\OO) := \inf\{t>0: X_t \notin \OO \}$.
Then define
\begin{equation}\label{e:2.3}
H_{\gamma} \varphi (x):=\< \gamma  , \nabla  \tilde{u}  \> (x),  \qquad x\in \po
. %\gamma_i (x) \frac{\partial \tilde{u} (x)}{ \partial x_i}.
\end{equation}

Denote by $ ( X^1, K^1)$ the solution of \eqref{e:2.1} with $ \vp=1$.
Since the local time $ K^1$ increases when and only when $ X^1$ hits the boundary $\po$, we can obtain a Markov process on the boundary by putting $ \widetilde{X}^1 (s):= X^1( K^{-1} (s) )$, where $ K^{-1} (s) $ is the right continuous inverse $ \sup\{t: K^1 (t) \leq s  \}$ of $ K^1$.
In \cite[Theorem 4]{Ueno}, %  divergence , time change
it shows that the operator $ H_{\gamma}$ is the generator of the process $\widetilde{X}^1 $.
By the periodicity of the coefficients, the Markov process $\widetilde{X}^1 $ induces a Markov process $\widetilde{X}_{\T^{d-1}}^1 $ on the torus $ \T^{d-1}$.
Combined with %In view of
the compactness of $ \T^{d-1}$ and Doeblin's theorem , we can deduce by a similar argument in \cite[Lemma 4.3]{Tanaka} that there exists a unique invariant measure $ \tilde{m}$ of the Markov process $\widetilde{X}_{\T^{d-1}}^1 $.
\bigskip

We now list some general assumptions for the semilinear PDEs \eqref{e:1.1}. %All these assumptions are assumed to hold in the sequel unless otherwise specified.

\begin{assumption}\label{h:2}
The functions $ a, b, c$ are all periodic of period 1 in each component. The coefficient $b$ is bounded and $ c$ satisfies
\begin{equation}\label{e:2.4}
-\infty <- \alpha \leq c(x) \leq 0 , \qquad \forall x\in \R^d,
\end{equation}
for some positive constant $\alpha $.
Moreover, we assume that for every $ 1\leq i \leq d$,
$ \sum_{j=1}^d  \frac{\partial a_{ij} }  { \partial x_j} (x) \in L^\infty (\R^d)$ and the following condition holds
\begin{equation}\label{e:2.5}
-\frac12 \sum_{j=1}^d  \int_{\T^d} a_{ij} (x) \, \frac{ \partial m (x)}{ \partial x_j}  \, dx +  \int_{\T^d}
b_i (x) m(x) \, dx =0.
\end{equation}
\end{assumption}

\begin{remark}\label{R:2.1} \rm
The condition \eqref{e:2.5} is common and natural in the homogenization problem and the following comments will be helpful to understand it. Note that $ a_{ij}$ is smooth enough, then

{\rm (i)} the operators $ L^\vp$ can be rewritten as
$$
L^\vp  = \frac12   \sum_{i,j=1}^d    a_{ij} ( x / \varepsilon)  \frac{\partial^2 }{ \partial x_i \partial x_j}  +
\frac{1}{\vp}    \sum_{i=1}^d  \left(
\frac12  \sum_{j=1}^d \frac{ \partial a_{ij}  }{  \partial x_j}  ( x/\vp )
+  b_i (x /\vp)
\right)\, \frac{ \partial }{ \partial x_i}.
$$
It is easy to see that \eqref{e:2.5} is the centering condition in the book of Bensoussan et al. \cite{Bensoussan} (see also \cite{B Hu, Lejay b, Pardoux, Tanaka}).
%and one can also find it in \cite{B Hu, Lejay b, Pardoux, Tanaka}.

{\rm (ii)} the operator $ L$ can be written as
$$
  L= \frac12  \sum_{i,j=1}^d   a_{ij}(x)     \frac{\partial^2}{ \partial x_i \partial x_j}  + \sum_{i=1}^d  \left(\frac12 \sum_{j=1}^d  \frac{\partial a_{ij} }  { \partial x_j} (x) +b_i(x) \right)  \, \frac{\partial}{ \partial x_i}.
$$
As %Following a similar argument
in \cite[Remark 4]{Tanaka}, \eqref{e:2.5} ensures each component process $X_i (t) $ of the $L$-diffusion $X(t)$ is recurrent in the sense that $X_i (t) $ hits any state in $\R^d$ with probability 1.

{\rm (iii)} in the Neumann boundary case, \eqref{e:2.5} corresponds to the condition (H.3) in \cite{Tanaka}, and in the Dirichlet boundary case corresponds to Eq.(38) in \cite{Lejay b}.
%Let $ \wh b_i (x):= (\frac12 \sum_{j=1}^d  \frac{\partial a_{ij} }  { \partial x_j}  +b_i) (x)$, then
Moreover, under the condition \eqref{e:2.5}, there will be a unique periodic solution to the auxiliary problem \eqref{e:3.1} in Section \ref{S:3} for each $ i=1,\dots, d$.
Hence the solution of the problem \eqref{e:3.1} can be given by $ \omega_i(x)= \int_0^\infty E_x [ \wt b_i ( \overline{X}_t)] \, dt$ when $ \overline{X}$ starts from $x \in \T^d$ at time $0$.
\end{remark}

%In this paper, we always extend the definition of the functions $f, c$ and $g$ to $\R^d$ by setting their values to be zero off $\OO$.

\begin{assumption}\label{h:3}
The function $ f: \R^d \times \R \times \R^d \to \R$ is a bounded uniformly continuous function which satisfies
\begin{description}
  \item{\rm (1)} $(y_1-y_2)( f(x,y_1,z)-  f(x,y_2,z)  ) \leq c_1(x) \, |y_1-y_2|^2 $, where $ c_1(x) $ is a Borel measurable function.

  \item{\rm (2)} $|f(x,y,z_1)- f(x,y,z_2)  | \leq c_2 \, |z_1-z_2 |$. Here $c_2$ is a positive constant.
\end{description}
\end{assumption}
Moreover, we assume %We also assume the following integrability condition hold
\begin{equation}\label{e:2.6}
E_x^\vp \left[   e^{ 2 \int_0^T   c_1^+ ( X^\vp(s) ) \, ds %+ 2( \alpha +1)\, K_T^\vp
}   \right]  < \infty,   % \qquad ~\hbox{for any fixed}~\vp,
\end{equation}
for any fixed $ \vp$ %with the non-negative constant $ \alpha$, %(see \eqref{e:2.4}).
%where $ \alpha$ is the non-negative constant in \eqref{e:2.4}.
%This condition will be used in the proof of Theorem \ref{T:1.1} in Section \ref{S:5}.
and also
\begin{equation}\label{e:2.7}
 E_{ x }^0 \Big[ e^{ 2\int_0^T c_1^+ (X^0 (r)) dr } \Big]  < \infty,
\end{equation}
holds where $E_x^0  $ denotes the expectation under the law of the reflected Brownian motion $ X^0$ with covariance matrix $ \bar {a}$. Here %Here we point out
\eqref{e:2.6} and \eqref{e:2.7} are imposed to ensure that the inequalities \eqref{e:5.17} and \eqref{e:5.18} in the proof of homogenization hold, respectively. In particular, they are clearly true when $ c_1$ is a negative or bounded function.

\begin{assumption}\label{h:4}
$g: \boo \to \R$ is continuous and bounded. In this paper, we always extend the definition of the function $g$ to $\R^d$ by setting its values to be zero off $\OO$.%\qquad or $ g\in L^2 (\po )$. }
\end{assumption}

\section{The auxiliary periodic problems}\label{S:3}
%we deal with the periodic solutions of the auxiliary problems \eqref{e:3.1}. Meanwhile, the $ L^p $-estimates $( p>2)$ of the gradient of the solution are obtained using the results of Gr\"{o}ger \cite{Groger} (see Theorem \ref{T:3.2}).
In this section, we study the periodic solutions of the auxiliary problems \eqref{e:3.1}. This will make us to get rid of the highly oscillating terms in treating the reflecting diffusion process $X^\vp $ in Section \ref{S:5}.
The main result is Theorem \ref{T:3.2}, which gives the $ L^p $-estimates $( p>2)$ of the gradient of the solution of the auxiliary problem.
We are mainly inspired by the thoughts of \cite[Theorem 1]{Groger} and %a similar argument as that in
\cite[Theorem 2]{Meyer}.

Based on \eqref{e:1.3} and Assumptions \ref{h:2}, we now
consider the 1-periodic solutions in $ H^1(\T^d)$ of the auxiliary problems
\begin{equation}\label{e:3.1}
\begin{cases}
L \omega_i =  - \left(\frac12 \sum_{j=1}^d  \frac{\partial a_{ij} }  { \partial x_j}  +b_i \right)
 %&\hbox{in the weak sense} ,
 \cr
\int_{\T^d} \omega_i(x) \, m(x) dx= 0 ,
\end{cases}
\end{equation}
in the weak sense for $i=1,\dots , d $. %This makes us to get rid of the highly oscillating terms in treating the reflecting diffusion process $X^\vp $ in Section \ref{S:5}.
Combined with \eqref{e:2.5} and Fredholm alternative theorem, it follows from \cite[Theorem 3.3.5]{Bensoussan} that the solution to \eqref{e:3.1} exists and is unique.
On the other hand, as stated in the proof of \cite[Proposition 1]{Lejay a} or \cite[Theorem 7.2]{S}, each function $ \omega_i$ is continuous and bounded. Define
\begin{eqnarray*}
\wt \omega_i (x) := x_i + \omega_i^\sharp (x),  \quad~\hbox{and}\quad
\wt \omega_i^\vp (x) := \vp \,\wt \omega_i (x /\vp),
\end{eqnarray*}
for $ i=1, \dots , d$, where $ \omega_i^\sharp$ is the extension of $\omega_i $ by periodicity to the whole $\R^d$.
Then in the weak sense, $\wt \omega_i $ satisfies %the functions $\wt \omega_i $ and $ \wt \omega_i^\vp$ satisfy
$$
\frac12 \nabla \cdot( A \, \nabla \wt \omega_i) +b \cdot \nabla \wt \omega_i =0.
$$
Meanwhile,
$$
\frac12 \nabla \cdot \big( A(x/\vp) \,\nabla \wt \omega_i^\vp (x)\big) +
\frac{1}{\vp}\, b(x/\vp) \cdot \nabla \wt \omega_i^\vp (x) =0,
$$
holds for each function $ \wt \omega_i^\vp$.
Hence $\wt \omega_i $ and $ \wt \omega_i^\vp$ are harmonic functions for operators $ L$ and $L^\vp$, respectively.
Since $\wt \omega^\vp_i $ belongs to the domain of the quadratic form associated with the process $ X^\vp$,
it follows from the Fukushima's decomposition (see \cite{FOT}) that
\begin{eqnarray}
d\wt \omega^\vp_i  ( X_s^\vp )
&=& \< \nabla \wt \omega_i  ( X_s^\vp / \vp)  , dM_s^\vp  \>
+ \< \nabla \wt \omega_i,  \gamma \> ( X_s^\vp / \vp)  \, dK_s^\vp            \nonumber \\
\label{e:3.2}
&=:& d \wt M_s^{\vp,i} + \wt \gamma_i   ( X_s^\vp / \vp)  \, dK_s^\vp ,
\qquad  \quad i=1, \dots , d, \quad 0 \leq t\leq s\leq T,
\end{eqnarray}
staring from $x$ at time $t$. %and we assume that $ \wt \omega^\vp_i  ( X_s^\vp )$ for any $ \vp >0$.
Moreover, $\wt M_s^{\vp} =( \wt M_s^{\vp,1}, \dots, \wt M_s^{\vp,d} )^T $ is a local martingale with cross-variations
\begin{eqnarray}
 \ll \wt M^{\vp,i}, \wt M^{\vp,j} \gg_s  &=&    \int_0^t
\< A \nabla \wt \omega_i,   \nabla \wt \omega_j \>  ( X_s^\vp / \vp)   \,ds  \nonumber \\
\label{e:3.3}
&=:& \int_0^t  \wh a_{ij} ( X_s^\vp / \vp)   \,ds.  %\nonumber \\
%&\xrightarrow{\vp \to 0}&  \left( \int_{\T^d}   \< A \nabla \wt \omega_i,   \nabla \wt \omega_j \>  ( x) \, m(x)dx
 %\right)  \, ds \nonumber \\
%&=:& \bar{a}_{ij} \, ds.
\end{eqnarray}
We take $ \nabla \wt \omega_i (x)$ as the column vectors to form a matrix, and denote it by  $\nabla \wt \omega (x)$. Let
$  \wt{\gamma}^\vp  (x):= (A    \nabla \wt \omega)^* (x /\vp) \, \vec{n} (x)$ %with $ \nabla \wt \omega (x)= (\nabla \wt \omega_1 (x), \dots,   \nabla \wt \omega_d (x))$.  Denote
and $ \wt \omega^\vp (x) := ( \wt \omega_1^\vp (x), \dots, \wt \omega_d^\vp (x) )^T$,
%$ \wt \omega_i^\vp $
%satisfies Eq.\eqref{e:3.2} and $ \wt \gamma_i (x/\vp) =  \< \gamma, \nabla \wt \omega_i \> (x/\vp) $ for every $ i=1, \dots , d$.
then it yields from \eqref{e:3.2} that
\begin{equation}\label{e:w}
d \wt \omega^\vp (X_s^\vp) = d \wt M_s^{\vp}+  \wt{\gamma}^\vp (X_s^\vp  )\, dK_s^\vp.
\end{equation}

\begin{remark}\label{R:3.1}
The operator $L= \frac12 \nabla \cdot( A\nabla) +b \cdot \nabla $ inside $\OO$ equipped with the Neumann boundary condition $ \<A(x)\vec{n} , \nabla \cdot \> =0$ on $\po$ determines the reflecting diffusion process $X^1$.
Then the scaling relation shows that $ \vp X^1 ( \cdot / \vp^2)$ is equivalent in law to $ X^\vp (\cdot)$. In addition, the measure $ m(x)dx$ is invariant for the reflected process $ X^1$. That is, we have
\begin{equation}\label{e:3:f}
\int_{\boo} E_x [ f(X^1 (t))]  \, m(x)dx
= \int_{\boo} f(x) \, m(x)dx,
\end{equation}
for any bounded and continuous function $f$ over $\boo$.

It suffices to prove \eqref{e:3:f} holds for the function $ f\in C_0^\infty (\OO)$.
Indeed, it is well known that
$ v(t,x):= E_x [ f(X^1 (t))]$ gives the probabilistic representation to the problem
\begin{eqnarray*}
\begin{cases}
\partial_t   v(t,x)= Lv(t,x), \quad
&  [0,t] \times \OO,  %x\in \po, ~t\in [0,T),
\cr
\frac{\partial v}{ \partial \gamma} (t,x)=0,  \quad
&  [0,t] \times \po,  %x\in \po, ~t\in [0,T),
\cr
v(0,x)= f(x), \quad
& x\in \OO.
\end{cases}
\end{eqnarray*}
then we obtain
\begin{eqnarray*}
\partial_t   \left(    \int_{\boo} v(t,x)  \, m(x)dx   \right)
&=& \int_{\boo}  Lv(t,x)  \, m(x)dx  \\
&=& \int_{\boo}  v(t,x)  \, L^* m(x)dx \\
&=&0.
\end{eqnarray*}
This implies $ \int_{\boo} v(t,x)  \, m(x)dx$ is constant with respect to $t$. In view of $ v(0,x)= f(x)$, it yields the equality \eqref{e:3:f} holds.

Based on the above displays, it should be pointed out that we can also use the thoughts of \cite[Proposition 2.4]{Pardoux} to establish the following convergence result
$$ \int_0^t  \wh a_{ij} ( X_s^\vp / \vp)  ds
\stackrel{dist}{=}
\vp^2  \int_0^{t/ \vp^2}  \wh a_{ij} ( X_s^1 )  ds
\xrightarrow{\vp \to 0}
t \int_{\T^d}   \< A \nabla \wt \omega_i,   \nabla \wt \omega_j \>  ( x) \, m(x)dx,
$$
%(homogenization stable lemma4.6 +  Bensoussan theorem 3.3.1 + pardoux 1999 proposition 2.1-2.4 + reflecting heat \\
%or reflecting heat   proposition 2.1 )
by applying the properties of $m$. More generally, we can prove the conclusion \eqref{e:5.7}, which plays an important role in homogenization in Section \ref{S:5}.
\end{remark}

%Define $  \bar{\Gamma} := ( \bar{\Gamma}_1,
%\dots, \bar{\Gamma}_d )$ with
%\begin{equation}\label{e:3.4}
%\color{red}
%\bar{\Gamma}_i = \int_{\T^{d-1}  }  \< \gamma, \nabla \wt \omega_i \> (x)  \, \wt m (dx)   \qquad i=1,\dots, d.
%\end{equation}
%Combined with the above results and \cite[Theorem 2.1]{Tanaka}, then we can obtain the following theorem which plays an important role in our study of the homogenized equation.
%\begin{thm}\label{T:3.1}
%Under the Assumptions \ref{h:2} and \eqref{e:1.3}, the $ (L^\vp ,  \partial_{\gamma_\vp})$-diffusion $ X^\vp$  converges in law to a $ (\overline{L}, \partial_{ \bar{\Gamma}  } ) $-diffusion $ \wh X$ as $ \vp \downarrow 0$. Moreover,
%$$
%(   X^\vp, \wt M^\vp, K^\vp) \Rightarrow ( \wh X, \wh M , \wh K ),
%$$
%where $\wh  M, \wh K$ are the martingale part and local time of $\wh X$, respectively.
%\end{thm}

\begin{thm}\label{T:3.2}
Let $ B(0, R)$ be the ball of $\R^d$ centered at $0$, of radius $R>0$. Define the function $\mathcal{G} : B(0, 2R) \times \R^{d+1} \to    \R^{d+1} $,
\begin{eqnarray*}
\mathcal{G}_j (x, \zeta)&:=&  \frac12 \sum_{i=1}^d   a_{ij}(x)  \zeta_i,  \qquad\quad j=1, \cdots , d,  \\
\mathcal{G}_0 (x, \zeta)&:=& \sum_{i=1}^d   b_i(x)  \zeta_i,
\end{eqnarray*}
for $ \zeta= ( \zeta_0, \dots , \zeta_d)^T \in \R^{d+1}$
and assume the following conditions
\begin{equation}
\begin{array} {lllllc}
\<\mathcal{G} (x, \zeta) -  \mathcal{G} (x,  \vartheta ), \zeta-  \vartheta   \> \,  \geq  q_1 |\zeta-\vartheta|^2, \quad q_1>1/(2\lambda),
\\
\label{e:3.5}  \\
|  \mathcal{G} (x, \zeta) -  \mathcal{G} (x,  \vartheta ) |  \leq   q_2 |\zeta-\vartheta|, \quad q_2<
\infty,
%\quad \hbox{for}~ x\in B(0, 2R),  \zeta,\vartheta \in \R^{d+1}.
\end{array}
\end{equation}
hold for any $x\in B(0, 2R) $ and $ \zeta,\vartheta \in \R^{d+1}$.
%Assume the conditions
%\begin{equation}
%\begin{array} {lllllc}
%\<b, \nabla \zeta- \nabla \vartheta   \> \, (\zeta-\vartheta)  \geq \lambda^{-1} |\zeta-\vartheta|^2
%\\
%\label{e:3.G}  \\
%|  \<b, \nabla \zeta- \nabla \vartheta   \>  |  \leq \lambda |\zeta-\vartheta|,
%\end{array}
%\end{equation}
%hold for any $\zeta,\vartheta \in H^1 (\T^d ) $.
Then there exists a constant $Q(\lambda, d ,q_1, q_2) >2  $ such that for all $p\in [2, Q(\lambda, d ,q_1, q_2) ) $, $\nabla \omega_i \in L^p ( \T^d) ,~i=1, \cdots d$.
\end{thm}

\pf
%Let $ B(0, R)$ be the ball of $\R^d$ centered at $0$, of radius $R>0$.
Setting
$\wt   \omega_{i, R} (x):= \wt  \omega_i(x) - \int_{ B(0, 2R) }  \wt \omega_i  (x) dx $, then it is easy to see that
$ \int_{ B(0, 2R) }  \wt  \omega_{i, R} (x) dx =0 $ and each function $ \wt \omega_{i, R}   \in H_{loc}^1 ( \R^d)$ satisfies
$$
\frac12 \nabla \cdot( A \, \nabla \wt \omega_{i, R}) +b \cdot \nabla \wt \omega_{i, R} =0 ,  \qquad i=1, \cdots d,
$$
in the weak sense.
In view of \eqref{e:1.3} and \eqref{e:3.5}, it yields that the function $\mathcal{G} %: B(0, 2R) \times \R \times \R^d \to    \R^{d+1}
$
%defined by
%\begin{eqnarray*}
%\mathcal{G}_j (x, \zeta, \nabla \zeta)&:=& \sum_{i=1}^d   a_{ij}(x) \frac{ \partial \zeta}{\partial x_i} ,  \qquad\quad j=1, \cdots , d,  \\
%\mathcal{G}_0 (x, \zeta, \nabla \zeta)&:=& \sum_{i=1}^d   b_i(x) \frac{ \partial \zeta}{\partial x_i},
%\end{eqnarray*}
satisfies the conditions (4.1) of \cite{Groger}.
For $ u\in H^1 (B(0, 2R) )$, define the operator $\Lambda\in {\mathcal{L} (  H^1; L^2 ( B(0, 2R); \R^{d+1}) )  }$ by $ \Lambda u:= ( u, \nabla u)^T$.
Let the operator $ J: W_0^{1,2} (B(0, 2R)) \to W^{-1,2} (B(0, 2R) )$ be %defined by
\begin{eqnarray}
\forall  v \in  W_0^{1,2}, \qquad
\<Ju,v \> &:=& \int_{ B(0, 2R)}  \mathcal{G}(\cdot, \Lambda u) \cdot \Lambda v \, dx \nonumber \\
\label{e:3.7}
&=& \int_{ B(0, 2R)}  \bigg( \frac12 \sum_{i,j=1}^d a_{ij} (x) \frac{\partial u}{ \partial x_i}  \frac{\partial v}{ \partial x_j}
+ \sum_{i=1}^d  b_i(x)  \frac{\partial u}{ \partial x_i}  v(x)
\bigg)\, dx.
\end{eqnarray}
%is strongly monotone and Lipschitzian.
Hence by \cite[Theorem 1]{Groger} and a similar argument as that in \cite[Theorem 2]{Meyer}, we can deduce the existence of constants $ Q(\lambda, d ,q_1, q_2)>2$ and $ C( \lambda, p, d)>0$ such that for all $ p\in [2, Q(\lambda, d ,q_1, q_2) )$,
\begin{equation}\label{e:3.6}
\big\|  \nabla  \wt \omega_{i, R}    \big\|_{ p, B(0,R) }
\leq C( \lambda, p, d)   \,  R^{ d  \big(\frac1p- \frac12 \big) -1}  \,  \big\| \wt \omega_{i, R}    \big\|_{2, B(0,2R)},
\end{equation}
where $ \| \cdot \|_{ p, B(0,R) } $ denotes the norm in $ L^p ( B(0,R), dx)$.

On the other hand, in view of $ \int_{ B(0, 2R) }  \wt  \omega_{i, R} (x) dx =0 $,
it follows from Poincar\'{e} inequality that
%$$
%\big\| \wt \omega_{i, R}     -
%\int_{ B(0, 2R) }  \wt  \omega_{i, R} (x) dx
%\big\|_{2, B(0,2R)}   \leq C(d) \, R \, \big\|  \nabla  \wt \omega_{i, R}    \big\|_{ 2, B(0,2R) },
%$$

$$
\big\| \wt \omega_{i, R}    \big\|_{2, B(0,2R)}
\leq C(d) \, R \, \big\|  \nabla  \wt \omega_{i, R}    \big\|_{ 2, B(0,2R) }.
$$
Hence combined with \eqref{e:3.6}, we have
\begin{eqnarray*}
R^{ - \frac{d}{p} }    \big\|  \nabla  \wt \omega_{i, R}    \big\|_{ p, B(0,R) }    \leq
C( \lambda, p, d)   \,  R^{ - \frac{d}{2} } \big\|  \nabla  \wt \omega_{i, R}    \big\|_{ 2, B(0,2R) }   .
\end{eqnarray*}
Moreover, since
\begin{eqnarray*}
&&\limsup_{R\to \infty}  R^{ - d }   \,  \big\|  \nabla  \wt \omega_{i, R}    \big\|_{ 2, B(0,2R) }^2  \\
&=&\limsup_{R\to \infty}   \frac{1}{R^d}  \int_{ B(0, 2R) }   \big\| \nabla \wt \omega_{i, R} (x)    \big\|^2   \, dx     \\
&=& \limsup_{R\to \infty}    \frac{1}{R^d}  \int_{ B(0, 2R) }   \big\|    e_i  + \nabla  \omega_{i}^\sharp (x)    \big\|^2   \, dx     \\
&=&  |B(0,2)|    \int_{\T^d}   \big\| e_i  + \nabla  \omega_{i} (x)    \big\|^2   \,  dx,
\end{eqnarray*}
and
\begin{eqnarray*}
&&\limsup_{R\to \infty}   R^{ - d }  \,  \big\|  \nabla  \wt \omega_{i, R}    \big\|_{ p, B(0,R) }^p   \\
&=&\limsup_{R\to \infty}   \frac{1}{R^d}  \int_{ B(0, R) }   \big\| \nabla \wt \omega_{i, R} (x)    \big\|^p   \, dx     \\
&\geq& \limsup_{R\to \infty}    \frac{1}{R^d}  \int_{ B(0, R) }   \big[  \big\|    e_i  + \nabla  \omega_{i}^\sharp (x)    \big\|  \wedge n \big]^p   \, dx     \\
&=&  |B(0,1)|    \int_{\T^d}   \big[ \big\| e_i  + \nabla  \omega_{i} (x)    \big\|    \wedge n \big]^p   \,  dx,
\end{eqnarray*}
%we can obtain a conclusion from monotone convergence theorem that
then monotone convergence theorem implies
$$
\big\| e_i  + \nabla  \omega_{i} (x)    \big\|_{L^p (\T^d)}
\leq  C( \lambda, p, d) \big\| e_i  + \nabla  \omega_{i} (x)    \big\|_{L^2 (\T^d)} < \infty.
$$
The proof is complete.
\qed

\begin{remark}\label{R:3.2}
The conditions \eqref{e:3.5} mean that the operator $J$ defined by \eqref{e:3.7} is strongly monotone and Lipschitzian.
More precisely,
$ ( \zeta_0 - \vartheta_0) ( \sum_{i=1}^d  b_i(x)  (\zeta_i - \vartheta_i )  ) \geq (q_1- 1/(2\lambda) ) \, |\zeta - \vartheta|^2$ and
$ \big( \sum_{i=1}^d  b_i(x)  (\zeta_i - \vartheta_i )  \big)^2  + \sum_{j=1}^d  \big( \sum_{i=1}^d   a_{ij}(x) (\zeta_i - \vartheta_i )
  \big)^2  \leq   5 q_2^2 |\zeta - \vartheta|^2$.
Clearly, \eqref{e:3.5} holds when the coefficients $ a_{ij}$ and $b$ are bounded.
%In view of \eqref{e:1.3}, then \eqref{e:3.5} is equivalent to the following conditions hold
%for any $\zeta,\vartheta \in H^1 (\T^d ) $. That is to say
\end{remark}

\section{Backward SDEs with singular coefficients}\label{S:4}
%In this section, we mainly study the existence and uniqueness of solutions for a class of BSDEs, which involves the integral with respect to the local time of a reflecting diffusion process.    In Section \ref{S:4}, we mainly study the existence and uniqueness of a class of BSDEs with singular coefficients (see Theorem \ref{T:4.1}).
This section is independent of other sections and devoted to study the existence and uniqueness of a class of BSDEs with singular coefficients, which involves the integral with respect to the local time of a reflecting diffusion process. The main result is Theorem \ref{T:4.1} and it implies in Section \ref{S:5} that for any fixed $\vp$, there exists a unique pair $ (Y_s^\vp,  Z_s^\vp )_{ s\in [ t,T]}$ of progressively measurable processes satisfying \eqref{e:5.4} and \eqref{e:5.5}.

For any fixed $ \vp$, let $ ( \Omega, \P, \FF_t)$ be the probability space carrying the reflecting diffusion process $ X(t), t\geq 0$ described in Section \ref{S:2} and $ M(t), K(t)$ are the martingale part and local time of $X$, resectively.
By the martingale representation theorem in \cite[Theorem 2.1]{Zhang}, we mainly study the existence and uniqueness of solutions for a class of BSDEs associated with the martingale part $M(t)$ and the local time $ K(t) $.

Let $ F(\omega, s,y,z):\Omega \times [0,T] \times \R \times \R^d  \to \R$ and $ h(\omega, s,y): \Omega \times [0,T] \times \R \to \R$ be given progressively measurable functions. For simplicity, we omit the random parameter $\omega$. Assume that they are continuous in $y$ and satisfy the following conditions
\begin{description}
  \item{\rm (A.1)} $(y_1-y_2)( F(s,y_1,z)-  F(s,y_2,z)  ) \leq d_1(s) |y_1-y_2|^2 $,

  \item{\rm (A.2)} $(y_1-y_2)( h(s,y_1)-  h(s,y_2)  ) \leq
   \beta(s) |y_1-y_2|^2 $,

  \item{\rm (A.3)} $|F(s,y,z_1)- F(s,y,z_2)  | \leq d_2 |z_1-z_2 |
  $,

  \item{\rm (A.4)} $  | F(  s,y, z) | \leq |F(s,0,z)| + d_3(s) (1+|y|) $,

  \item{\rm (A.5)} $| h(s, y) | \leq |h(s,0)|+d_3(s) ( 1+|y| )$,
\end{description}
where $ d_1(s) , d_3(s) $ are progressively measurable stochastic process, $d_2$ is a positive constant and $ %- \bar{\beta} \leq
\beta(s) < 0 $ for all $s\in [0,T]$. %for some constant $ \beta $.
Let $ \xi \in L^2(\Omega, \FF_T, P)$.

\begin{lemma}\label{L:4.1}
Denote
$$
 \varphi (s):= \int_0^s  d(u) du+ \int_0^s \mu (u)  \,dK_u,
$$
%denote $ \int_0^t  d(u) du+ \int_0^t \mu (u) dK_u$ by $ \varphi (t)$
where $d(s):= 2d_1^+(s)  $ and $\mu (s):= 2 ( \beta(s)+1) $.
Let
$ E\big[ e^{ \varphi(T)}   |\xi |^2  \big] < \infty$, $ E\big[  \int_0^T   e^{ \varphi(s) }  \,  |h(s,0)|^2 dK_s  \big]   < \infty$ and
$$
E\left[ \int_0^T    e^{ \varphi(s)}  \, \Big(  |F(s,0,0)|^2  + | d_3 (s)|^2  \Big) ds
\right] <\infty,
$$
then there exists a unique solution $ ( Y,Z)$ to the following BSDE
\begin{equation}\label{e:4.1}
Y(t) = \xi + \int_t^T  F(s, Y(s), Z(s)) \, ds + \int_t^T  h(s,Y(s)) \, dK_s  - \int_t^T  \<Z(s) , dM(s)\>.
\end{equation}
Furthermore, %the solution $  ( Y,Z)$  satisfies
\begin{equation}\label{e:4.12}
E\Big[ \sup_{t \in [0,T] }   e^{ \varphi (t)}   |Y(t)|^2    \Big]  < \infty,
\qquad \qquad
E\Big[   \int_0^T   e^{ \varphi (s)}
\| Z(s)\|^2\,ds \Big]  < \infty,
\end{equation}
and
\begin{equation}\label{e:4.13}
E\Big[   \int_0^T   e^{ \varphi (s)}     |Y(s) |^2    \,  dK_s
\Big]  < \infty.
\end{equation}
\end{lemma}

\pf
{\bf Uniqueness:}
%Define $ d(s) := -2 d_1(s)$ and
Suppose $ ( Y^1(t), Z^1 (t)), ( Y^2(t), Z^2 (t))$ are two solutions to Eq.\eqref{e:4.1}.
%Then
%\begin{eqnarray*}
%d( | Y^1(t)- Y^2(t) |^2 )
%&=& -2 (Y^1(t)- Y^2(t)  )  \, \big(   F( t,Y^1(t), Z^1(t) )  -F( t,Y^2(t), Z^2(t) )   \big) \, dt  \\
%&&- 2 (Y^1(t)- Y^2(t)  ) \,  \big(   h( t,Y^1(t) )  -h( t,Y^2(t) )   \big) \, dK_t  \\
%&& +  2 (Y^1(t)- Y^2(t)  )  \,   \< Z^1(t)- Z^2(t) , dM(t)  \>  \\
%&&+ \big\< \wh a( X_t) ( Z^1(t)- Z^2(t)  ),  Z^1(t)- Z^2(t)  \big\> \, dt
%\end{eqnarray*}
By Ito's formula and (A.1)-(A.3), we have
\begin{eqnarray*}
&&e^{ \varphi (t) %\int_0^t d (s) ds
}  |  Y^1(t)- Y^2(t) |^2
+ \int_t^T   e^{ \varphi (s) }   \big\<   a( X(s))  ( Z^1(s)- Z^2(s)  ),  Z^1(s)- Z^2(s)    \big\>  \, ds
\\
&=&  -  \int_t^T  d(s)  \, e^{ \varphi (s) }   |  Y^1(s)- Y^2(s) |^2  \, ds   - \int_t^T  \mu(s)  \, e^{ \varphi (s) }   |  Y^1(s)- Y^2(s) |^2  \, dK_s  \\
&&+2\int_t^T e^{ \varphi (s) }  (   Y^1(s)- Y^2(s)  ) \,
\big( F(s, Y^1(s), Z^1(s) )  -F(s, Y^2(s), Z^2(s))  \big) \,ds
\\
&&+2\int_t^T e^{ \varphi (s) }  (   Y^1(s)- Y^2(s)  ) \,
\big( h(s, Y^1(s) )  -h(s, Y^2(s))  \big) \, dK_s
\\
&&- 2\int_t^T e^{ \varphi (s) } (   Y^1(s)- Y^2(s)  ) \,
\< Z^1(s)- Z^2(s), dM(s)   \>
\\
&\leq&  -  \int_t^T  d(s)  \, e^{ \varphi (s) }   |  Y^1(s)- Y^2(s) |^2  \, ds   + 2\int_t^T   d_1(s)\,  e^{ \varphi (s) }\,
|   Y^1(s)- Y^2(s) |^2 ds  \\
&&+  2\int_t^T   d_2 \, e^{ \varphi (s) }\,  | Y^1(s)- Y^2(s) |     | Z^1(s)- Z^2(s)  |  ds - 2 \int_t^T   e^{ \varphi (s) }   |  Y^1(s)- Y^2(s) |^2  \, dK_s
\\
&&- 2\int_t^T e^{ \varphi (s) } (   Y^1(s)- Y^2(s)  ) \,
\< Z^1(s)- Z^2(s), dM(s)   \>
\\
&\leq&    \frac12   \int_t^T  \frac{1}{\lambda}  e^{ \varphi (s) }  \,
|  Z^1(s)- Z^2(s) |^2 ds
+ 8d_2^2 \lambda    \int_t^T e^{ \varphi (s) }  \,  |   Y^1(s)- Y^2(s) |^2 ds  \\
&& - 2\int_t^T e^{ \varphi (s)} (   Y^1(s)- Y^2(s)  ) \,
\< Z^1(s)- Z^2(s), dM(s)   \>
 -2 \int_t^T    e^{ \varphi (s) }   |  Y^1(s)- Y^2(s) |^2  \, dK_s .
\end{eqnarray*}
%In view of (A.1) and (A.3), it yields
%\begin{eqnarray*}
%&&2\int_t^T e^{ \int_0^s d(u) du }  (   Y^1(s)- Y^2(s)  ) \,
%\big( F(s, Y^1(s), Z^1(s) )  -F(s, Y^2(s), Z^2(s))  \big) ds
%\\
%&=& 2\int_t^T e^{ \int_0^s d(u) du }  (   Y^1(s)- Y^2(s)  ) \,
%\big( F(s, Y^1(s), Z^1(s) )  -F(s, Y^2(s), Z^1(s))  \big) ds
%\\
%&&+ 2\int_t^T e^{ \int_0^s d(u) du }  (   Y^1(s)- Y^2(s)  ) \,
%\big( F(s, Y^2(s), Z^1(s) )  -F(s, Y^2(s), Z^2(s))  \big) ds
%\\
%&\leq&  - 2\int_t^T   d_1(s)\,  e^{ \int_0^s d(u) du }\,
%|   Y^1(s)- Y^2(s) |^2 ds  \\
%&&+  2\int_t^T e^{ \int_0^s d(u) du }\, d_2 \, | Y^1(s)- Y^2(s) |     | Z^1(s)- Z^2(s)  |  ds
%\\
%&\leq&   - 2\int_t^T   d_1(s)\,  e^{ \int_0^s d(u) du }\,
%|   Y^1(s)- Y^2(s) |^2 ds
%+  \frac12   \int_t^T e^{ \int_0^s d(u) du }  \, \lambda
%|  Z^1(s)- Z^2(s) |^2 ds
%\\
%&&+ \frac{8d_2^2}{ \lambda}    \int_t^T e^{ \int_0^s d(u) du }  \,  |   Y^1(s)- Y^2(s) |^2 ds ,
%\end{eqnarray*}
%where we have used Young's inequality in the last step.
Thus combined with \eqref{e:1.3} and (A.2), we can obtain
\begin{eqnarray}
&&e^{ \varphi (t)}  |  Y^1(t)- Y^2(t) |^2
+ \frac12  \int_t^T   e^{ \varphi (s) }   \big\<   a( X(s))  ( Z^1(s)- Z^2(s)  ),  Z^1(s)- Z^2(s)    \big\>  \, ds
\nonumber   \\
\label{e:4.2}
&\leq&  -  2\int_t^T e^{ \varphi (s) } (   Y^1(s)- Y^2(s)  ) \,   \< Z^1(s)- Z^2(s), dM(s)   \>
\nonumber   \\
&& + 8d_2^2 \lambda    \int_t^T e^{ \varphi (s) }  \,  |   Y^1(s)- Y^2(s) |^2 ds .
\end{eqnarray}
Taking expectation in the above inequality, it yields
$$
E \left[  e^{\varphi (t) }  |  Y^1(t)- Y^2(t) |^2
\right] \leq
C_\lambda    \int_t^T  E\left[   e^{ \varphi (s) }  \,  |   Y^1(s)- Y^2(s) |^2     \right]  \, ds.
$$
By Gronwall's inequality, we conclude that $  \forall t,  Y^1(t)=Y^2(t), \, a.s.$ and hence $ Z^1(t)= Z^2(t), \, a.s.$ by \eqref{e:4.2}.
\bigskip

{\bf Existence:} Define
$
F_n (t,y,z):= \int_{\R} F(t,x,z) \phi_n (y-x ) dx
$ and $
h_n ( t,y):= \int_{\R} h(t,x) \phi_n (y-x ) dx,
$
where $ \phi_n (x):=n \phi(nx)$ and $ \phi \in C_0^\infty (\R) $ is an even, nonnegative function with $ \int_{\R} \phi(x) dx=1$.
Hence, it is easy to see that for each $n\geq 1$, %$ F_n$ and $h_n$ satisfy  %Lipschitz in $y$, uniformly with respect to $(t,)$
\begin{eqnarray}\label{e:4.3}
|F_n (t,y_1,z)- F_n (t,y_2,z)| \leq C_n |y_1- y_2|, \quad
|h_n(t,y_1)- h_n(t,y_2)| \leq C'_n |y_1- y_2|,
\quad y_1,y_2 \in \R,
\end{eqnarray}
for some constants $C_n $ and $ C'_n$. Furthermore, since functions $ F$ and $h$ are continuous in $y$, we have $F_n (t,y,z) \to F (t,y,z) $ and $h_n (t,y) \to h (t,y) $ as $n \to \infty$.

Consider the following BSDE
\begin{eqnarray}\label{e:4.4}
Y_n(t) = \xi + \int_t^T F_n (s,Y_n(s), Z_n(s)   ) ds +
\int_t^T h_n (s,Y_n(s)  ) dK_s - \int_t^T \<Z_n (s), dM(s)\>.
\end{eqnarray}
In view of \eqref{e:4.3} and the assumptions (A.3)-(A.5), we deduce from \cite[Theorem 1.6]{Pardoux Zh} that \eqref{e:4.4} admits a unique solution $ (Y_n, Z_n)$.
Now, our aim is to show that there exists a convergent subsequence $ ( Y_{n_k}, Z_{n_k}  )$.
Indeed by Ito's formula, it yields
\begin{eqnarray*}
&&e^{ \varphi (t)}   |Y_n(t)|^2
+ \int_t^T    e^{ \varphi (s)}
\mu (s)   |Y_n(s)|^2 \, dK_s     + \int_t^T   e^{ \varphi (s)}
\<  a(X(s)) Z_n (s), Z_n(s)     \> \,ds
\\
&=&   e^{ \varphi (T)}   |\xi |^2
- \int_t^T   e^{\varphi (s)}  \,  d(s) |Y_n(s)|^2 \, ds
+ 2 \int_t^T   e^{ \varphi (s)}  \,   Y_n(s)  \,   F_n (s, Y_n(s),  Z_n(s) )\, ds
\\
&&+ 2 \int_t^T   e^{ \varphi (s)}  \,   Y_n(s)  \,  h_n (s, Y_n(s) ) \,  dK_s
- 2  \int_t^T   e^{ \varphi (s)}  \,   Y_n(s)  \,   \<Z_n(s), dM(s)    \>.
\end{eqnarray*}
By (A.1-4), \eqref{e:1.3} and Young's inequality, we have
\begin{eqnarray*}
&&2 \int_t^T   e^{ \varphi (s)}  \,   Y_n(s)    F_n (s, Y_n(s),  Z_n(s) )\, ds      \\
%&\leq&  -2  \int_t^T   e^{ \int_0^s  d(u) du+ \int_0^s \mu (u) dK_u}  \, d_1(s)  | Y_n(s)  |^2   \,ds
%\\
%&&+  2 d_2  \int_t^T   e^{ \int_0^s  d(u) du+ \int_0^s \mu (u) dK_u}  \,   | Y_n(s)  |    | Z_n(s)  | ds
%\\
%&&+  2   \int_t^T   e^{ \int_0^s  d(u) du+ \int_0^s \mu (u) dK_u}  \,   | Y_n(s)|  | F(s,0,0) | ds
%\\
&\leq&    2  \int_t^T   e^{ \varphi (s)}  \, d_1(s)  | Y_n(s)  |^2   \,ds   + \frac12 \int_t^T   e^{ \varphi (s)}
\<  a(X(s)) Z_n (s), Z_n(s)     \> \,ds
\\
&&+ \Big( 8 d_2^2\lambda+1  \Big)   \int_t^T   e^{ \varphi (s)}  \,   | Y_n(s)  |^2    ds
+ \int_t^T   e^{ \varphi (s)}  \,   | F(s,0,0)  |^2    ds ,
\end{eqnarray*}
and
\begin{eqnarray*}
&&2 \int_t^T   e^{ \varphi (s)}  \,   Y_n(s)    h_n (s, Y_n(s) ) \,  dK_s
\\
&\leq& 2  \int_t^T   e^{ \varphi (s)}  \,   \beta(s) |Y_n(s) |^2    \,  dK_s
+ 2 \int_t^T   e^{ \varphi (s)}  \,   |Y_n(s)|    |h (s,0 )|^2  \,  dK_s
\\
&\leq& \int_t^T   e^{\varphi (s)}  \,    ( 2 \beta(s) +1)   |Y_n(s) |^2    \,  dK_s
+ \int_t^T   e^{ \varphi (s)}  \, |h(s,0)|^2  \,  dK_s .
\end{eqnarray*}
Consequently,
\begin{eqnarray}
&&  e^{ \varphi (t)}   |Y_n(t)|^2
+ \int_t^T    e^{ \varphi (s)}   |Y_n(s)|^2 \, dK_s
+ \frac12  \int_t^T   e^{ \varphi (s)}
\<  a(X(s)) Z_n (s), Z_n(s)     \> \,ds
\nonumber  \\  \label{e:4.5}
&\leq& e^{ \varphi (T)}   |\xi |^2
+C_{ \lambda}  \int_t^T   e^{ \varphi (s)}  \,   | Y_n(s)  |^2    ds     + \int_t^T   e^{ \varphi (s)}  \,   | F(s,0,0)  |^2    ds  \nonumber \\
&& + \int_t^T   e^{ \varphi (s)}  \, |h(s,0)|^2  \,  dK_s
- 2  \int_t^T   e^{ \varphi (s)}  \,   Y_n(s) \<Z_n(s), dM(s)    \>.
\end{eqnarray}
Taking expectation and by Gronwall's inequality, it yields
\begin{eqnarray}
&&\sup_n \sup_{0\leq t \leq T}  E\left[ e^{ \varphi (t)}   |Y_n(t)|^2  \right]
\nonumber \\    \label{e:4.6}
&\leq& C\bigg\{  E\left[ e^{\varphi (T)}   |\xi |^2   \right]
+   E\left[ \int_t^T   e^{ \varphi (s)}  \,   | F(s,0,0)  |^2    ds  \right]     +E\left[ \int_t^T   e^{ \varphi (s)}  \, |h(s,0)|^2  \,  dK_s   \right]
\bigg\}     \nonumber \\
&<& \infty.
\end{eqnarray}
Hence, we deduce that
\begin{equation}\label{e:4.7}
\sup_n E\left[   \int_t^T   e^{ \varphi (s)}
\<  a(X(s)) Z_n (s), Z_n(s)     \> \,ds
\right]  < \infty,
\end{equation}
and
\begin{equation}\label{e:4.8}
\sup_n E\left[   \int_t^T   e^{ \varphi (s)}     |Y_n(s) |^2    \,  dK_s
\right]  < \infty.
\end{equation}

Furthermore, combined with the conditions of this lemma, we also obtain from \eqref{e:4.5} and \eqref{e:4.6} that there exists some constant $C$ such that
\begin{eqnarray*}
&& E\Big[ \sup_{t \in [0,T] }   e^{ \varphi (t)}   |Y_n(t)|^2    \Big]
\nonumber  \\     \label{e:4.9}
&\leq&   E\left[  e^{ \varphi (T)}   |\xi |^2  \right]
+C  E\left[  \int_0^T   e^{ \varphi (s)}  \,   | Y_n(s)  |^2    ds  \right]
\nonumber \\
&& + E\left[  \int_0^T  e^{ \varphi (s)}  \,   | F(s,0,0)  |^2    ds  \right]
+ E\left[  \int_0^T  e^{\varphi (s)}  \, |h(s,0)|^2  \,  dK_s   \right]
\nonumber \\
&&+C  E\Big[ \sup_{t\in [0,T]}  \int_0^t   e^{ \varphi (s)}  \,   Y_n(s) \<Z_n(s), dM(s)    \>
\Big]   \nonumber \\
&\leq& C+ CE\Big[ \sup_{t\in [0,T]}  \int_0^t   e^{ \varphi (s)}  \,   Y_n(s) \<Z_n(s), dM(s)    \>    \Big]  .
\end{eqnarray*}
By both Burkh\"{o}lder's and Young's inequalities, it implies
\begin{eqnarray*}
&&CE\bigg[ \sup_{t\in [0,T]}  \int_0^t   e^{ \varphi (s)}  \,   Y_n(s) \<Z_n(s), dM(s)    \>    \bigg]  \\
&\leq&  CE\bigg[ \bigg( \int_0^T   e^{ 2 \varphi (s)}  \,   Y_n^2(s)
\big\<  a(X(s)) Z_n (s), Z_n(s)     \big\> \,ds  \bigg)^{\frac12}  \bigg]   \\
%&\leq&  CE\bigg[  \bigg(  \sup_{s\in [0,T]}    e^{ \frac12 \int_0^s  d(u) du+ \frac12 \int_0^s \mu (u) dK_u}  \,  | Y_n(s)|   \bigg)
%\\
%&&\qquad ~\times   \bigg(  \int_0^T   e^{ \int_0^s  d(u) du+ \int_0^s \mu (u) dK_u}  \,    \big\< \wh a(X(s)) Z_n (s), Z_n(s)     \big\> \,ds     \bigg)^{\frac12}
%\bigg]     \\
&\leq&  \frac12   E\Big[ \sup_{s \in [0,T] }   e^{ \varphi (s)}  \, |Y_n(s)|^2    \Big]
+  CE\bigg[   \int_0^T   e^{ \varphi (s)}  \,    \big\< a(X(s)) Z_n (s), Z_n(s)     \big\> \,ds  \bigg].
\end{eqnarray*}
Then it follows from \eqref{e:4.7} that
\begin{equation}\label{e:4.10}
\sup_n  E\Big[ \sup_{t \in [0,T] }   e^{ \varphi (t)}   |Y_n(t)|^2    \Big]  < \infty.
\end{equation}
In view of \eqref{e:4.7}-\eqref{e:4.10}, %\eqref{e:4.8} and \eqref{e:4.10},
we can extract a subsequence $ n_k$ such that $ Y_{n_k}(t) e^{ (1/2)\varphi (t)} $ converges to some $ \wh Y(t)$ in $L^2 ( \Omega, L^\infty[0,T])$ equipped with the weak star topology. In addition, $ Z_{n_k}(t) e^{ (1/2)\varphi (t)}$ converges weakly to some $ \wh Z(t)$ in $L^2 ( [0,T]\times \Omega; \R)$.
Since
\begin{eqnarray*}
e^{ \frac12  \varphi (t)}   Y_{n_k}(t)
&=& e^{ \frac12  \varphi (T)}  \xi +  \int_t^T   e^{ \frac12  \varphi (s)}  \,  F_{n_k} (s,Y_{n_k}(s),  Z_{n_k}(s)  )  \, ds
\\
&& + \int_t^T   e^{ \frac12  \varphi (s)} \,   h_{n_k} (s,Y_{n_k}(s)  )  \, dK_s
- \frac12  \int_t^T   e^{ \frac12  \varphi (s)}  \,   Y_{n_k}(s)  d(s) \, ds
\\
&& - \frac12  \int_t^T   e^{ \frac12 \varphi (s)}  \,   Y_{n_k}(s)  \mu(s) \, dK_s
- \int_t^T   e^{ \frac12  \varphi (s)} \,    \<Z_{n_k}(s), dM(s)    \>,
\end{eqnarray*}
letting $k\to \infty$ and by the same arguments in the proof of \cite[Proposition 1.8]{Pardoux Zh} (see P.546-547), we conclude that the limit $ ( \wh Y, \wh Z  )$ satisfies
\begin{eqnarray*}
\wh   Y(t)
&=&  e^{ \frac12 \varphi (T)}  \xi   +  \int_t^T   e^{ \frac12  \varphi (s)}
\times  \,  F \big(s,  e^{ -\frac12  \varphi (s)}   \wh Y(s),
 e^{ -\frac12  \varphi (s)}  \wh Z(s)  \big)  \, ds
\\
&& + \int_t^T   e^{ \frac12  \varphi (s)} \,   h (s, e^{ -\frac12  \varphi (s)}   \wh Y(s)  )  \, dK_s
 - \int_t^T     \<\wh Z(s), dM(s)    \>
\\
&&- \frac12  \int_t^T  \wh Y(s)    d(s) \, ds
 - \frac12  \int_t^T   \wh Y(s)  \mu(s) \, dK_s
.
\end{eqnarray*}
Define
$
Y(t) := e^{ -(1/2)  \varphi (t)} \,  \wh Y( t),
$ and $
Z(t) := e^{ -(1/2) \varphi (t)} \,  \wh Z( t),
$
then Ito's formula yields that
$$
Y(t)= \xi + \int_t^T  f (s, Y(s), Z(s) ) ds
+  \int_t^T  h (s, Y(s) ) dK_s   - \int_t^T    \<\wh Z(s), dM(s)    \>,
$$
which implies $(Y,Z)$ is a solution to the backward equation \eqref{e:4.1}. Applying Fatou's lemma, \eqref{e:4.12} and \eqref{e:4.13} follows from the above proof.
%The proof is complete.
\qed

Now we apply Lemma \ref{L:4.1} to a particular situation.
Let $ F(x,y,z): \R^d \times \R \times \R^d  \to \R$ and $ h(x,y): \R^d \times \R$ be Borel measurable functions. Assume that they are continuous in $y$ and satisfy the following conditions
\begin{description}
  \item{\rm (B.1)} $(y_1-y_2)( F(x,y_1,z)-  F(x,y_2,z)  ) \leq d_1(x) |y_1-y_2|^2 $,

  \item{\rm (B.2)} $(y_1-y_2)( h(x,y_1)-  h(x,y_2)  ) \leq
   \beta(x) |y_1-y_2|^2 $,

  \item{\rm (B.3)} $|F(x,y,z_1)- F(x,y,z_2)  | \leq d_2 |z_1-z_2 |
  $,

  \item{\rm (B.4)} $  | F(  x,y, z) | \leq |F(x,0,z)| + d_3(x) (1+|y|) $,

  \item{\rm (B.5)} $| h(x, y) | \leq |h(x,0)|+d_3(x) ( 1+|y| )$,
\end{description}
where $ d_1$ and $ d_3$ are Borel measurable functions on $\R^d$, $d_2$ is a positive constant and $ \beta$ is a bounded negative measurable function on $\R^d$.
Given $ g \in C_b (\R^d)$ and consider the BSDE
\begin{equation}\label{e:4.11}
Y(t) = g(X(T)  ) + \int_t^T  F(X(s) , Y(s), Z(s)) \, ds + \int_t^T  h(X(s),Y(s)) \, dK_s  - \int_t^T  \<Z(s) , dM(s)\>,
\end{equation}
where $ M(s)$ is the martingale part of $ X(s)$. Set $ d(x):= 2d_1^+(x) $ and $ \mu (x):= 2( \beta (x)+1 )$, the following results follows from Lemma \ref{L:4.1}.
\begin{thm}\label{T:4.1}
Let (B.1)-(B.5) hold. Assume moreover
$
E\big[ e^{ \int_0^T d( X(s) ) ds+ \int_0^T \mu ( X(s) ) dK_s}   \big] < \infty
$,
$$ E\left[  \int_0^T   e^{ \int_0^s d( X(u) ) du  + \int_0^s \mu ( X(u)  ) dK_u}  \,  |h(X(s) ,0)|^2 dK_s  \right]   < \infty,
$$
and
$$
E\left[ \int_0^T    e^{ \int_0^s d( X(u) ) du  + \int_0^s \mu ( X(u) ) dK_u}  \, \Big(  |F(X(s) ,0,0)|^2  + | d_3 (X(s) )|^2  \Big) ds
\right] <\infty,
$$
then the BSDE \eqref{e:4.11} admits a unique solution.
\end{thm}

\section{Homogenization of parabolic systems}\label{S:5}
In this section, we are concerned with the homogenization of the parabolic systems \eqref{e:1.1}. In Subsection \ref{S:5.1}, the homogenized coefficients are defined (see \eqref{e:a}-\eqref{e:5.2}) and we show that the homogenized boundary value problem \eqref{e:5.3} has a unique weak solution (Theorem \ref{T:5.1}).
In Subsection \ref{S:5.2}, inspired by the thoughts of \cite[Lemma 6.1, Lemma 6.3]{Tanaka}, two important lemmas are proved (see Lemma \ref{L:5.2} and Lemma \ref{L:5.3}).
Based on them, the homogenization result can be obtained in the case where the coefficients $ f $ and $ g$ are smooth (Lemma \ref{L:5.4}). At the end of Subsection \ref{S:5.3}, Theorem \ref{T:1.1} is proved by a regularization procedure.

\subsection{Homogenized PDEs with the third boundary conditions}\label{S:5.1}
%Let us define the coefficients
Define
\begin{equation}\label{e:a}
 \bar{a}_{ij}:=
\int_{\T^d}   \< A \nabla \wt \omega_i,   \nabla \wt \omega_j \>  ( \eta) \, m(\eta)  d\eta ,
\end{equation}
\begin{eqnarray}
\, \bar{f} (x,y,z)&:=& \int_{\T^d} f\big(x,y, \nabla\wt \omega ( \eta)  \, z  \big)\, m(\eta) d\eta
\nonumber   \\
\label{e:5.1}
&=& \int_{\T^d} f\big(x,y, (Id+  \nabla \omega ( \eta)  ) \, z  \big)\, m(\eta) d\eta,
\end{eqnarray}
and
\begin{equation}\label{e:5.2}
\bar{C}:= \int_{\T^{d-1}  } c(\eta) \, \wt m(d\eta).
\end{equation}
\begin{thm}\label{T:5.1}
%Assume \eqref{e:1.3} and \ref{h:2}-\ref{h:4} are satisfied.
Let \eqref{e:1.3} and Assumptions \ref{h:2}-\ref{h:4} hold,
then $ \bar{A} $ defined by \eqref{e:a} is a strictly positive symmetric matrix, and $\bar f$ satisfies Assumption \ref{h:3} with constant $ c_2  \int_{\T^d} \| (Id+  \nabla \omega ( \eta) \|   \, d\eta  $.
Moreover, the homogenized boundary value problem
\begin{eqnarray}\label{e:5.3}
\begin{cases}
\frac{ \partial u^0}{ \partial t} (t,x)+ \overline{L} u^0 (t,x) +  \bar{f}(x,u^0 (t,x),\nabla u^0 (t,x) )=0   & (t,x) \in [0,T) \times \OO,  %x\in\OO, ~t\in [0,T),
\cr
\frac12 \frac{\partial u^0}{ \partial  \upsilon^0    }  (t,x)
+\bar{C} u^0 (t,x)=0   & (t,x) \in [0,T) \times \po,  %x\in \po, ~t\in [0,T),
\cr
u^0 (T,x)= g(x) & x\in \boo,
\end{cases}
\end{eqnarray}
has a unique weak solution in the space
\begin{eqnarray*}
\mathcal{W}_1^2  (0,T, H^1(\OO),   L^2(\OO) )
&:=& \big\{  u(t,x)\in L^2([0,T]; H^1(\OO)  )   ~\hbox{such that}~ \\
&& \qquad \partial_t u(t,x)\in L^2([0,T], H^{-1}(\OO)  )
\big\}.
\end{eqnarray*}
\end{thm}

\pf
In view of $ \bar{a}_{ij}=\int_{\T^d}   \< A \nabla \wt \omega_i,   \nabla \wt \omega_j \>  ( x) \, m(x)dx$, $\bar{A}$ is clearly a non-negative symmetric matrix. Hence, it suffices to prove that $\bar{A}$ is non-degenerate. Assume that there exists $ \xi \in \R^d$ such that $ \< \bar{A}\, \xi, \xi\> =0$. Using the fact that $ A$ is uniformly elliptic, we can have $\int_{\T^d}  \<  \nabla \wt \omega_i (x), \xi\> \, dx =0$ for every $ i=1 , \cdots, d$. That is
$$
\xi_i + \sum_{ j=1}^d   \int_{\T^d}  \big( \nabla  \omega_i \big)_j  (x)   \,  \xi_j \, dx =0 ,
\qquad \forall   i=1 , \cdots, d.
$$
In view of $ \omega \in  H^1(\T^d)$, then it yields $  \int_{\T^d}  \nabla  \omega_i (x)  \, dx=0$.  This implies that
$\xi=0$.

The assertion concerning $\bar f$ is an easy consequence of Assumption \ref{h:3} and $\nabla \omega_i \in L^2 ( \T^d) $ for $i=1, \cdots d$. Since the coefficient $\bar A$ is a constant matrix, the existence and uniqueness of a weak solution to \eqref{e:5.3} can be deduced from \cite[Theorem 1]{WYZ}.
\qed

\subsection{Two lemmas}\label{S:5.2}
Under Assumptions \ref{h:2}-\ref{h:4}, one can deduce by a similar argument as that in \cite[Proposition 3.2]{Pardoux Zh} together with Khas'minskii's lemma that for any fixed $\vp$,
$
E_x^\vp [ \,  \exp \{  2( \alpha +1)\, K_T^\vp
\}   ]  < C(\alpha, T).
$
%In view of \eqref{e:1.3} and
Hence by \eqref{e:2.6}, Theorem \ref{T:4.1} and the boundedness of function $f$,
%it shows for any fixed $\vp$ that
there exists a unique pair $ (Y_s^\vp,  Z_s^\vp )_{ s\in [ t,T]}$ of progressively measurable processes satisfying
%solution $ (Y_s^\vp,  Z_s^\vp )_{ s\in [ t,T]}$ to the following backward SDEs
\begin{eqnarray}
Y^\vp (s) &=& g (X^\vp (T) ) + \int_s^T  f( X^\vp (r), Y^\vp (r), Z^\vp (r)  ) \, dr   + \int_s^T  c( X^\vp (r) / \vp )  Y^\vp (r) \, dK_r^\vp     \nonumber \\
\label{e:5.4}
&& - \int_s^T  \< Z^\vp (r) , dM^\vp (r)  \>, \quad t\leq s \leq T , \qquad  \P_{t,x}^{\vp} -a.s. ,
\end{eqnarray}
and
\begin{equation}\label{e:5.5}
E_{ t,x }^\vp  \left[  \sup_{t\leq s \leq T}   | Y^\vp (s)  |^2   + \int_t^T  \| Z^\vp (s)  \|^2 ds
\right] < \infty.
\end{equation}
Moreover, for each fixed $\vp$, we deduce from \cite[Corollary 3]{WYZ} that $ Y^\vp (s)= u^\vp (  s,X^\vp (s) )$ where $ u^\vp$ is a continuous version of the weak solution of system \eqref{e:1.1}. Therefore $ u^\vp (t,x) = Y^\vp (t )  $.
Let $ E_x^\vp$ be the expectation under $ \P_{t,x}^{\vp}$, we are going to prove that for all $ p\in \big( 1,  Q(\lambda, d ,q_1, q_2)/ 2   \big)$,
$$ \lim_{\vp \to 0} E_x^\vp [ |  Y^\vp (t )  - u^0  (t,x) | ^p ] =0,
$$
where $ Q(\lambda, d ,q_1, q_2)   $ is the constant in Theorem \ref{T:3.2}.
To avoid heavy notations, we will take in all the sequel $t=0$. To this end, we prove the following two lemmas, which will play an important role in the homogenization of Theorem \ref{T:5.2}.

\begin{lemma}\label{L:5.2}
Let $\psi(s,x ,\eta) : [0,T] \times \R^d \times \R^d \to \R$ be a bounded, uniformly continuous function which is periodic in $ \eta$. Assume that for any $ x\in \po$,
\begin{equation}\label{e:n}
%\color{blue}
\Big\< \Big( \int_{ \T^{d-1} }  ( A \nabla \wt \omega)_{ij} ( \eta)   \, \wt m (d\eta)  \Big) \,\nabla \Psi (x),   \nabla \Psi (x)     \Big\>
> 0,   %\qquad \forall   x\in \po,
\end{equation}
where the matrix $ (\int_{ \T^{d-1} }  ( A \nabla \wt \omega)_{ij} ( \eta)   \, \wt m (d\eta)\, )_{1\leq i,j \leq d }$ is constant
and matrix $\nabla \wt \omega $ is composed of $ \nabla \wt \omega_i$ as column vectors, that is $ \nabla \wt \omega (x)=  (\nabla \wt \omega_1 (x), \nabla \wt \omega_2 (x),  \cdots, \nabla \wt \omega_d (x))$. If
$
\int_{ \T^d}  \psi (s,x, \eta) \, m(\eta) d\eta =0
$
holds for any $ s\in [0,T]$ and $ x\in \R^d$,
then we have
\begin{equation}\label{e:5.6}
\lim_{\vp \to 0} \,
E_{ t,x }^\vp  \left[   \left|  \int_0^s
\psi \bigg(r, X_r^\vp , \frac{X_r^\vp}{ \vp }  \bigg) \, dr
\right|^2     \right] =0.
\end{equation}
\end{lemma}

\pf
By a standard smooth approximation procedure, it is sufficient to prove \eqref{e:5.6} holds when the function $ \psi$ is $ C_b^\infty$ and periodic in $ \eta$.
Indeed, we consider the equation
\begin{eqnarray*}
\begin{cases}
\frac12     \sum_{i,j=1}^d
\frac{\partial}{ \partial \eta_i}
\left(a_{ij} ( \eta)  \frac{\partial \varphi}{ \partial \eta_j} \right)  +     \sum_{i=1}^d
b_i (\eta)\frac{\partial  \varphi}{ \partial \eta_i}
= - \psi,
 \cr
\int_{\T^d} \varphi(s,x, \eta) \, m(\eta) d\eta= 0 .
\end{cases}
\end{eqnarray*}
By the regularity on $ \psi$ and on the coefficients, then we can use Ito's formula to obtain
\begin{eqnarray*}
\varphi \bigg(s, X_s^\vp , \frac{X_s^\vp}{ \vp }  \bigg)
&=& \varphi  \bigg(t, X_t^\vp , \frac{X_t^\vp}{ \vp }   \bigg)
+ \int_t^s  \frac{ \partial \varphi }{  \partial \tau}   \bigg( r,  X_r^\vp , \frac{X_r^\vp}{ \vp }   \bigg)  \, dr
+ \int_t^s    \bigg\< \bigg(  \frac{\partial \varphi }{ \partial x}   +\frac{1}{\vp}  \frac{\partial \varphi }{ \partial \eta}  \bigg) , d M_r^\vp     \bigg\>
\\
&&+ \int_t^s  \bigg( \frac{1}{\vp}  \Big\< \frac{\partial \varphi }{ \partial x}  , \wt b   \Big\>  +    \frac{1}{\vp^2} \, L\varphi   \bigg)  \, dr      + \int_t^s    \bigg\<  \bigg(  \frac{\partial \varphi }{ \partial x}   +\frac{1}{\vp}  \frac{\partial \varphi }{ \partial \eta}  \bigg) , \gamma(  X_r^\vp /\vp ) \bigg\> \, dK_r^\vp
\\
&&+ \frac12   \int_t^s  tr\bigg( a(  X_r^\vp /\vp )  \bigg(  \frac{\partial^2 \varphi }{ \partial x^2}  + \frac{2}{\vp}  \frac{\partial^2 \varphi }{ \partial x  \partial \eta}  \bigg) \bigg)\, dr.
\end{eqnarray*}
Hence
\begin{eqnarray*}
\int_t^s \psi \bigg(r, X_r^\vp , \frac{X_r^\vp}{ \vp }  \bigg) \,dr
&=&  \vp^2 \left[  \varphi  \bigg(t, X_t^\vp , \frac{X_t^\vp}{ \vp }   \bigg)  -  \varphi \bigg(s, X_s^\vp , \frac{X_s^\vp}{ \vp }  \bigg)     \right]   \\
&& + \vp^2 \int_t^s  \frac{ \partial \varphi }{  \partial \tau}   \bigg( r,  X_r^\vp , \frac{X_r^\vp}{ \vp }   \bigg)  \, dr
+ \int_t^s    \bigg\< \bigg( \vp^2 \frac{\partial \varphi }{ \partial x}   +\vp  \frac{\partial \varphi }{ \partial \eta}  \bigg) , d M_r^\vp     \bigg\>    \\
&&+ \vp  \int_t^s     \Big\< \frac{\partial \varphi }{ \partial x}  , \wt b   \Big\> \, dr   +    \int_t^s    \bigg\< \bigg( \vp^2 \frac{\partial \varphi }{ \partial x}   +\vp  \frac{\partial \varphi }{ \partial \eta}  \bigg)  , \gamma(  X_r^\vp /\vp ) \bigg\> \, dK_r^\vp    \\
&&+   \frac{\vp^2}{2}   \int_t^s  tr\bigg( a(  X_r^\vp /\vp )  \frac{\partial^2 \varphi }{ \partial x^2}  \bigg)\, dr
+\vp \int_t^s   tr\bigg( a(  X_r^\vp /\vp )     \frac{\partial^2 \varphi }{ \partial x  \partial \eta}  \bigg) \, dr.
\end{eqnarray*}

If $ E_{t,x}^\vp  [( K^\vp (T) )^2 ]  $ is bounded in $\vp$ for any fixed $ T>0$, it is easy to see that
$$
E_{ t,x }^\vp  \left[   \left|  \int_0^s
\psi \bigg(r, X_r^\vp , \frac{X_r^\vp}{ \vp }  \bigg) \, dr
\right|^2     \right] \leq C \vp,
$$
which implies \eqref{e:5.6} follows.
%On the other hand, based on the thoughts of \cite[Lemma 5.2, Lemma 6.1]{Tanaka}, we now prove that $ E_{t,x}^\vp  [( K^\vp (T) )^2 ]  $ is bounded in $\vp$ for any fixed $ T>0$.
In fact, based on the thoughts of \cite[Lemma 5.2, Lemma 6.1]{Tanaka}, we can obtain that $ E_{t,x}^\vp  [( K^\vp (T) )^p ]  $ is bounded in $\vp$ for each $ p\geq 1$.
For each function $\wt \omega_i$, consider the solution $ \phi_i$ of the following equation
\begin{eqnarray*}
\begin{cases}
\frac12   \sum_{i,j=1}^d
\frac{\partial}{ \partial x_i}
\left(a_{ij} ( x)  \frac{\partial  \phi}{ \partial x_j} \right)  +     \sum_{i=1}^d
b_i (x)\frac{\partial  \phi}{ \partial x_i}
= 0,   & x\in \OO,
 \cr
\frac{\partial \phi}{\partial \upsilon} (x)=\big \< A(x) \nabla \wt \omega_i (x)-
\big( \int_{ \T^{d-1} }  (A  \nabla \wt \omega_i)_j  ( \eta)   \, \wt m (d\eta)  \big)    , \vec{n}(x)  \big\>
& x\in \po ,
\end{cases}
\end{eqnarray*}
where the $\big( \int_{ \T^{d-1} }  (A  \nabla \wt \omega_i)_j  ( \eta)   \, \wt m (d\eta)  \big)_{1\leq j\leq d} $ is a constant vector.
As $ \phi_i \in W^{1,2} (\OO)$, it yields
\begin{eqnarray*}
d( \vp  \phi_i ( X_s^\vp  /\vp ) ) & = &  \< \nabla   \phi_i ( X_s^\vp /\vp ), dM_s^\vp\>  +   \<( A    \nabla \wt \omega_i) (X_s^\vp /\vp), \vec{n}(X_s^\vp )\>  \\
&&-
\Big\<  \Big( \int_{ \T^{d-1} }  (A  \nabla \wt \omega_i)_j  ( \eta)  \, \wt m (d\eta)  \Big), \vec{n}(X_s^\vp )   \Big\>.
\end{eqnarray*}
Let $ \phi (x):= (\phi_1(x), \cdots, \phi_d (x))^T$, then we have
\begin{eqnarray*}
\int_0^t   ( A \nabla \wt \omega)^* ( X_s^\vp  /\vp )  \, \vec{n}(X_s^\vp ) dK_s^\vp
&=&\int_0^t    \Big( \int_{ \T^{d-1} }  ( A \nabla \wt \omega)_{ij} ( \eta)   \, \wt m (d\eta)  \Big)^*   \, \vec{n}(X_s^\vp ) dK_s^\vp   \\
&& + \vp \big( \phi (  X_t^\vp  /\vp )   - \phi (x/\vp) \big)
- \int_0^t  \< \nabla   \phi ( X_s^\vp /\vp ), dM_s^\vp \>   .
\end{eqnarray*}
%where the matrix $ (\int_{ \T^{d-1} }  ( A \nabla \wt \omega)_{ij} ( \eta)   \, \wt m (d\eta)\, )_{1\leq i,j \leq d }$ is constant
%and matrix $\nabla \wt \omega $ is composed of $ \nabla \wt \omega_i$ as column vectors, that is
%$$ \nabla \wt \omega (x)=  (\nabla \wt \omega_1 (x), \nabla \wt \omega_2 (x),  \cdots, \nabla \wt \omega_d (x)) . $$
Consequently, combined with the definition of $ \wt \omega_i^\vp$ and \eqref{e:w}, it implies
\begin{eqnarray}\label{e:skorohod}
X^\vp_t  &=&  \wh X^\vp_t  +
\int_0^t    \Big( \int_{ \T^{d-1} }  ( A \nabla \wt \omega)_{ij} ( \eta)   \, \wt m (d\eta)  \Big)^*   \, \vec{n}(X_s^\vp ) dK_s^\vp,
\end{eqnarray}
where
\begin{eqnarray*}
\wh X_t^\vp &=& x + \big[   \wt M_t^\vp  -\int_0^t    \< \nabla   \phi ( X_s^\vp /\vp ), dM_s^\vp \>    \big]
 + \vp [\phi^\vp  (X_t^\vp  /\vp)  -\phi^\vp  (x/\vp) ]
\\
&&- \vp [ \wt \omega^\vp  (X_t^\vp  /\vp)  -\wt \omega^\vp  (x/\vp)  ]   .
\end{eqnarray*}
In view of \eqref{e:o} and \eqref{e:n}, then we can obtain %for any $ x\in \po$,
$$
\Big\< \Big( \int_{ \T^{d-1} }  ( A \nabla \wt \omega)_{ij} ( \eta)   \, \wt m (d\eta)  \Big) \,\vec{n} (x),   \vec{n} (x)     \Big\>
> 0,   \qquad \forall   x\in \po,
$$
which implies that \eqref{e:skorohod} can be regarded as a Skorohod equation with respect to $ \wh X^\vp (t)$ and $ K^\vp (t)$. From \cite[Theorem 2.2]{Costantini}, we can know that $ K^\vp (T)$ can be controlled by $  \sup_{0\leq t_1 \leq t_2 \leq T}    | \wh X^\vp (t_1) - \wh X^\vp (t_2) |$. Hence, $ E_{t,x}^\vp  [( K^\vp (T) )^p ]  $ is bounded in $\vp$ for any $ p\geq 1$.
The proof of this lemma is completed.
\qed

%{\color{blue}
%\begin{remark}\label{R:5.1}
%Here, we point out that the condition \eqref{e:n} is natural. For example,
%\end{remark}
%}

\begin{remark}\label{R:5.1}
 It should be noted that the method of \cite[Lemma 3.9.2]{Bensoussan} is no longer applicable because of the appearance of the local time term and the highly oscillating term $ \vp^{-1}  \wt b ( X^\vp /\vp) \, ds$ in \eqref{e:2.1}.
Under condition \eqref{e:n} and inspired by \cite[Lemma 6.1]{Tanaka}, process $X^\vp$ can be rewritten to another equation \eqref{e:skorohod} of Skorohod type. Hence from the property of the local time, we can prove $ E_{t,x}^\vp  [( K^\vp (T) )^p ]  $ is bounded in $\vp$ for any $ p\geq 1$.
The lemma is further proved to be true.

In addition, since $\wt \omega^\vp( X_s^\vp )= X_s^\vp + \vp ( \omega ( X_s^\vp  /\vp)- \omega(x/\vp)  ) $ and each function $ \omega_i$ is bounded, we can also conclude that %from Lemma \ref{L:5.2} that %the following corollary.
\begin{equation}\label{e:5.7}
E_{ t,x }^\vp  \left[   \left|  \int_0^s
\psi \bigg(r,\wt \omega^\vp( X_s^\vp ) , \frac{X_r^\vp}{ \vp }  \bigg) \, dr
\right|^2     \right]
\stackrel{\vp \to 0}{\longrightarrow}0.
\end{equation}
\end{remark}

Meanwhile, we have the following similar convergence result for integrals of the local time $K^\vp $.

\begin{lemma}\label{L:5.3}
Let $h : \R^{d} \to \R$ be bounded%or $ L^2$
, continuous and periodic of period one in each direction.
If the function $h$ satisfies $ \int_{\T^{d-1} } h (\eta) \, \wt m ( d\eta)=0$, then
\begin{equation}\label{e:5.9}
\lim_{\vp \to 0} \,
E_{ t,x }^\vp  \left[   \left|  \int_0^s
h \bigg( \frac{X_r^\vp}{ \vp }  \bigg) \, dK_r^\vp
%-\int_{\T^{d-1} } h (\eta) \, \wt m ( d\eta)
\right|^2     \right] =0.
\end{equation}
\end{lemma}

\pf
By a standard smooth approximation procedure, it suffices to prove \eqref{e:5.9} holds for the function $ h\in  C_b^\infty$. Let $ \phi$ be the solution of $L\phi=0 $ in $ \OO$ and $ \partial \phi/ \partial \gamma= h$ on $\po$.
%\begin{eqnarray*}
%\begin{cases}
%L\phi=0   &\hbox{in } \OO,   \cr
%\frac{\partial \phi}{\partial \gamma}= h & \hbox{on }\po,
%\end{cases}
%\end{eqnarray*}
Then by Ito's formula %and \eqref{e:2.1},
we have
\begin{eqnarray}
 \vp \phi (X_s^\vp/ \vp )
&=& \vp   \phi (x/ \vp )  + \frac{1}{\vp} \int_0^s  L\phi (X_r^\vp/ \vp ) \, dr   + \int_0^s  \<  \nabla \phi (X_r^\vp/ \vp ) , dM_r^\vp  \>  \nonumber \\
\label{e:5.10}
&&+ \int_0^s    \< \nabla \phi (X_r^\vp/ \vp ), \gamma (X_r^\vp/ \vp ) \>  \, dK_r^\vp   \nonumber \\
&=& \vp   \phi (x/ \vp )   +\int_0^s  \<  \nabla \phi (X_r^\vp/ \vp ) , dM_r^\vp  \>   + \int_0^s   h(X_r^\vp/ \vp ) \, dK_r^\vp
\nonumber \\
&=:&  \vp   \phi (x/ \vp )   +  I_{1,\vp} (s)  + I_{2,\vp} (s),
\end{eqnarray}
which implies
$$
E \left[  |  I_{1,\vp} (s)  + I_{2,\vp} (s)  |^2
\right]   = \vp^2 E[ |  \phi (X_s^\vp/ \vp )-  \phi (x/ \vp ) |^2 ]   \leq 4 \vp^2 | \phi|_{\infty}  <\infty.
$$

Now let $ Q_\vp$ be the probability measure induced by the %$ \R^2$-valued
process $( I_{1,\vp} (s)  , I_{2,\vp} (s)  )$, then as the arguments in \cite[Lemma 6.3]{Tanaka}, we can choose a subsequence such that $ Q_\vp $ converges weakly to some limit probability measure $ Q_0$ as $\vp \to 0$. Moreover, it follows from \eqref{e:5.10} that the limit $ (I_1(s), I_2 (s) )$
satisfies $ I_1(s) + I_2(s)=0$ ($Q_0$ a.s.) and they are a $ Q_0$-martingale and a bounded variation process, respectively. Hence $ I_1(s) =I_2(s)=0$ ($Q_0$ a.s.). We have proved the conclusion of the lemma.
\qed

\subsection{Homogenization}\label{S:5.3}
Now, we are going to prove Theorem \ref{T:1.1} in the case where the coefficients $ f $ and $ g$ are smooth. To this end, we introduce the following assumptions (C.1) and (C.2).
\begin{description}
  \item{\rm (C.1)} $f$ is bounded and $ f(x,0,0)\in L^2 ( \R^d )$. Moreover, $f$ is Lipschitz, that is
      $$
      | f(x,y,z)- f( x',y',z')| \leq C ( |x-x'|+ |y-y'| + |z-z'|), \quad \forall x,x'\in\R^d ,  \, y,y' \in\R,\, z,z'\in\R^d  .
      $$
\item{\rm (C.2)} $g : \R^d  \to \R$ is a $ C^2$ function.
\end{description}
Then we know from \cite[Theorem 7.4]{LSU} that the system \eqref{e:5.3} has a unique classical solution $u^0 \in  C^{1,2} ([0,T] \times \OO; \R )$.
%the solution $u^0$ of system \eqref{e:5.3} belongs to $ C^{1,2} ([0,T] \times \OO; \R )$

\begin{thm}\label{T:5.2}
Let Assumption \ref{h:2}, \eqref{e:1.3} and (C.1-2) hold. Define
\begin{eqnarray*}
\begin{cases}
\wt  Y^\vp (s):= Y^\vp (s) - u^0 ( s, \wt \omega^\vp (X_s^\vp) ),      \cr
\wt  Z^\vp (s):=Z^\vp (s) -     \nabla \wt \omega (  X_s^\vp  / \vp ) \,  \nabla u^0 ( s, \wt \omega^\vp (X_s^\vp) )     ,
\end{cases}
\end{eqnarray*}
then we have
\begin{eqnarray}
\wt  Y^\vp (s)
&=& g(X_T^\vp )  -  g( \wt \omega^\vp (  X_T^\vp )  )
+ \int_s^T  F^\vp ( r, \wt  Y_r^\vp, \wt  Z_r^\vp  ) \, dr
\nonumber \\
\label{e:5.12}
&&+ \int_s^T   \wh h^\vp ( r, \wt  Y_r^\vp) \, dK_s^\vp
- \int_s^T  \big\<  \wt  Z_r^\vp , dM_s^\vp  \big\>  ,
\qquad \P_{t,x}^{\vp} -a.s.,
\end{eqnarray}
where
\begin{eqnarray*}
F^\vp  (s,y,z ) &:=&  f\big( X_s^\vp, y+u^0( s, \wt \omega^\vp (X_s^\vp) ) , z+  \,  \nabla \wt \omega (  X_s^\vp  / \vp )\,  \nabla  u^0( s, \wt \omega^\vp (X_s^\vp) )
\big)   \\
&&- \bar{f} \big( \wt \omega^\vp (X_s^\vp) , u^0( s, \wt \omega^\vp (X_s^\vp) ),  \nabla  u^0( s, \wt \omega^\vp (X_s^\vp) )    \big)   \\
&&+ \frac12    \left(
\big( \wh a( X_s^\vp / \vp)   -\bar{a}    \big)_{ij}
\frac{ \partial^2  u^0}{\partial x_i  \partial x_j}       ( s, \wt \omega^\vp (X_s^\vp) )
\right),   \\
\wh h^\vp  (s, y) & :=&
2 c(X_s^\vp / \vp)  \, \big(y+   u^0( s, \wt \omega^\vp (X_s^\vp) )    \big)     -  2 \bar{C}  u^0( s, \wt \omega^\vp (X_s^\vp) )
\\
&&+\big\< \nabla  u^0( s, \wt \omega^\vp (X_s^\vp) )
,  \wt \gamma (X_s^\vp  / \vp )  -  \bar{A} \vec{n} ( \wt \omega^\vp (X_s^\vp)) \big\>.
\end{eqnarray*}
Moreover for all $t\in [0,T]$, we can deduce that as $ \vp $ tends to $0$, $  E_{ x }^\vp \big[  \big|  \int_0^t    F^\vp (s,0,0)
\, ds  \big|^2     \big] \to 0$ and $ E_{ x }^\vp  \big[  \big|  \int_0^t   \wh h^\vp  (s, 0)  \, dK_s^{\vp}
\big|^2    \big] \to 0$.
\end{thm}

\pf
By Ito's formula, we have
\begin{eqnarray*}
&&d\big( Y^\vp (s) - u^0 ( s, \wt \omega^\vp (X_s^\vp) )   \big)  \\
&=&  -\bigg[ f( X_s^\vp, Y_s^\vp, Z_s^\vp   )   -
\bar{f} (\wt \omega^\vp (X_s^\vp) , u^0( s, \wt \omega^\vp (X_s^\vp) ),  \nabla  u^0( s, \wt \omega^\vp (X_s^\vp) )    )
\\
&&+ \frac12   \Big(
\big( \wh a( X_s^\vp / \vp)   -\bar{a}    \big)_{ij}
\frac{ \partial^2  u^0}{\partial x_i  \partial x_j}     \Big)
  \bigg] ds
- \big[  2 c(  X_s^\vp  / \vp ) Y_s^\vp - 2 \bar{C}  u^0( s, \wt \omega^\vp (X_s^\vp) )
\\
&& +  \big\< \nabla  u^0( s, \wt \omega^\vp (X_s^\vp) ) , \,  \wt \gamma (X_s^\vp  / \vp )  -  \bar{A} \vec{n} ( \wt \omega^\vp (X_s^\vp))  \big\>
\big] \, dK_s^\vp
\\
&&+ \left\< Z_s^\vp -  \nabla \wt \omega (  X_s^\vp  / \vp )  \,  \nabla  u^0( s, \wt \omega^\vp (X_s^\vp) )   \,
   , \, dM_s^\vp  \right\> ,
\end{eqnarray*}
which implies \eqref{e:5.12}.
%Let
%$$
%\wt  Y^\vp (s):= Y^\vp (s) - u^0 ( s, \wt \omega^\vp (X_s^\vp) ), \quad
%\wt  Z^\vp (s):=Z^\vp (s) - ( \nabla \wt \omega (  X_s^\vp  / \vp ) )^*    \nabla u^0 ( s, \wt \omega^\vp (X_s^\vp) ) ,
%$$
%then it yields
Moreover,
\begin{eqnarray*}
F^\vp  (s,0,0 )&=&  \psi (s,   \wt \omega^\vp (X_s^\vp),   X_s^\vp /\vp )    +   \Big[
f\big( X_s^\vp, u^0( s, \wt \omega^\vp (X_s^\vp) ) ,   \nabla \wt \omega (  X_s^\vp  / \vp )  \,
\nabla  u^0( s, \wt \omega^\vp (X_s^\vp) )
\big)   \\
&&-   f \big( \wt \omega^\vp (X_s^\vp) , u^0( s, \wt \omega^\vp (X_s^\vp) ),    \nabla \wt \omega (  X_s^\vp  / \vp )  \, \nabla  u^0( s, \wt \omega^\vp (X_s^\vp) )       \big)
\Big]     \\
&&=:  \psi (s,   \wt \omega^\vp (X_s^\vp),   X_s^\vp /\vp )
+ I_1^\vp (s)  ,
%F_1^\vp (s)   +  F_2^\vp (s)  +F_3^\vp (s),
\end{eqnarray*}
with
\begin{eqnarray*}
\psi (s, x, \eta) &:=&   \Big[
f \big( x , u^0( s, x ),   \nabla \wt \omega (  \eta)  \,  \nabla  u^0( s, x )       \big)
-   \bar{f}   \big( x , u^0( s, x ),  \nabla  u^0( s, x)    \big)       \Big]    \\
&&+ \frac12    \Big(
\big( \wh a( \eta)   -\bar{a}    \big)_{ij}
\frac{ \partial^2  u^0}{\partial x_i  \partial x_j}  (s,x)
\Big),
\end{eqnarray*}
and
\begin{eqnarray*}
\wh h^\vp  (s, 0)  &=&
 \big\< \nabla  u^0( s, \wt \omega^\vp (X_s^\vp) )
,  \wt \gamma (X_s^\vp  / \vp  )  -  \bar{A} \vec{n} ( \wt \omega^\vp (X_s^\vp)) \big\>     \\
&&+ 2 [c(  X_s^\vp  / \vp  )   - \bar{C}  ]    \, u^0( s, \wt \omega^\vp (X_s^\vp) ) .
%&=:&  I_2^\vp (s) + 2 [c(  X_s^\vp  / \vp  )   - \bar{C}  ]    \, u^0( s, \wt \omega^\vp (X_s^\vp) ) .
\end{eqnarray*}

Since function $ f$ is Lipschitz with respect to $x$ and each solution $ \omega_i$ of the auxiliary problems is bounded, then we obtain
\begin{eqnarray*}
E_{ x }^\vp  \left[  \bigg|  \int_0^t    I_1^\vp (s)
\, ds  \bigg|^2     \right]
&\leq& C E_{ x }^\vp     \left[    \int_0^t  \big\| X_s^\vp    -   \wt \omega^\vp (X_s^\vp)    \big\|^2
\, ds    \right]   \\
&=& C \vp^2  \int_0^t   E_{ x }^\vp    [ \| \omega^\sharp ( X_s^\vp /\vp ) \|^2 ] \, ds\\
&\rightarrow&  0,
\end{eqnarray*}
as $\vp \to 0$.
By the definition of the homogenized coefficients, we obtain
$
\int_{\T^d} \psi (s,x, \eta)\, m(\eta) d\eta =0.
$
Then from Lemma \ref{L:5.2} and \eqref{e:5.7}, it yields
\begin{equation}\label{e:5.8}
\lim_{\vp \to 0}  E_{ x }^\vp    \left[  \bigg|  \int_0^t    F^\vp (s,0,0)
\, ds  \bigg|^2     \right]   %\stackrel{\vp \to 0}{\longrightarrow}
=0.
\end{equation}
On the other hand, by \cite[Proposition 8.5]{Oxford}, we know that $ a(X_s^\vp  / \vp  ) \nabla \wt \omega (X_s^\vp  / \vp  ) $ converges weakly to the homogenized matrix $ \bar{a}$ in $  L^2 ( \OO)^d$. Hence letting $\varepsilon \to 0$, we obtain
\begin{eqnarray*}
&&\big\< \nabla  u^0( s, \wt \omega^\vp (X_s^\vp) )
,  \wt \gamma (X_s^\vp  / \vp  )  -  \bar{A} \vec{n} ( \wt \omega^\vp (X_s^\vp)) \big\>       \\
&=&  \big\< \nabla  u^0( s, \wt \omega^\vp (X_s^\vp) ),
\big( ( A \nabla \wt \omega) (X_s^\vp  / \vp  ) -  \bar{A}     \big)^*    \,\vec{n}( X_s^\vp  )     \,
\big\>     \\
&&  +  \big\< \nabla  u^0( s, \wt \omega^\vp (X_s^\vp) ),    \,
\bar{A}   \big( \vec{n}( X_s^\vp  ) -   \vec{n} ( \wt \omega^\vp (X_s^\vp))     \big)
\big\>     \\
& \to & 0.
\end{eqnarray*}
Combined with Lemma \ref{L:5.3}, it implies
\begin{equation}\label{e:5.11}
\lim_{\vp \to 0}  E_{ x }^\vp    \left[  \bigg|  \int_0^t   \wh h^\vp  (s, 0)  \, dK_s^{\vp}
\bigg|^2    \right] =0  .
\end{equation}
The theorem is proved.
\qed

\begin{lemma}\label{L:5.4}
Let Assumptions \ref{h:2}-\ref{h:4}, \eqref{e:1.3} and (C.1-2) hold, then %for all $ p\in \big( 1, \frac{ Q(\lambda, d ,q_1, q_2)}{2}   \wedge 2   \big)$,
$$
%\lim_{\vp \to 0}  \,  | u^\vp (0,x ) - u^0 (0,x) |^p    =0.
\lim_{\vp \to 0}  \,u^\vp (0,x ) = u^0 (0,x).
$$
\end{lemma}

\pf
Define
\begin{eqnarray*}
\wh Y_t^\vp &:=& \wt Y_t^\vp  + \int_0^t  F^\vp (r,0,0)  \, dr
+ \int_0^t   \wh h^\vp  (r, 0)  \, dK_r^{\vp}    \\
&=&  \left[ g(X_T^\vp)  - g (\wt \omega^\vp (X_T^\vp ) ) \right]
+ \int_0^T   F^\vp (r,0,0)  \, dr  + \int_t^T \left[ F^\vp (r, \wt Y_r^\vp  , \wt Z_r^\vp)   - F^\vp (r,0,0)
\right]   dr    \\
&&- \int_t^T   \big\< \wt Z_r^\vp , dM_r^\vp   \big\>  +  \int_0^T   \wh h^\vp  (r, 0)  \, dK_r^{\vp}   +   \int_t^T \left[   \wh h^\vp  (r, \wt Y_r^\vp  )  -    \wh h^\vp  (r, 0)
\right]    dK_r^{\vp}.
\end{eqnarray*}
%To prove the lemma,
Then it suffices to show that for all $ p\in \big( 1, \frac{ Q(\lambda, d ,q_1, q_2)}{2}   \wedge 2   \big)$,
$ %\lim_{\vp \to 0}
  |\wh Y_0^\vp|^p %E_{ x }^\vp  \big[ |\wh Y_0^\vp|^p  \big]
$ tends to $0$ as $ \vp \to 0 $.

Note that for all $ p\in \big( 1, \frac{ Q(\lambda, d ,q_1, q_2)}{2}   \big)$, Theorem \ref{T:3.2} implies
$$
E_{ x }^\vp    \left[  \int_0^T  |  F^\vp (r,0,0)  |^p  \, dr
\right]  \leq C E_{ x }^\vp    \left[  \int_0^T  \Big(  1+ \|\nabla\omega^\sharp \|^{2p}     \Big)
%\big(  X_r^\vp  / \vp  \big)
\, dr    \right]   < \infty.
$$
Moreover, by H\"{o}lder's inequality and in view of the boundedness of functions $u^0, \nabla u^0,  c$, it yields
\begin{eqnarray*}
E_{ x }^\vp    \left[ \left( \int_0^T     \wh h^\vp  (r, 0)    \, dK_r^{\vp}   \right)^p \right]
&\leq& E_{ x }^\vp    \left[ \left( \int_0^T    | \wh h^\vp  (r, 0) |^p    \, dK_r^{\vp}   \right)    \left(  \int_0^T 1 \,  dK_r^{\vp}      \right)^{p-1}
\right]    \\
&\leq& C E_{ x }^\vp    \left[( K_T^{\vp})^{p-1}   \cdot \left( \int_0^T  \Big(  1+ \|\nabla\omega^\sharp \|^{p}     \Big)
%\big(  X_r^\vp  / \vp  \big)
\,dK_r^{\vp}   \right)  \right]   \\
&\leq& CE_{ x }^\vp    \left[( K_T^{\vp})^2    \right] .
%\\   &<& \infty.
\end{eqnarray*}
Hence combined with \eqref{e:5.5}, then for all $ p\in \big( 1, \frac{ Q(\lambda, d ,q_1, q_2)}{2}   \wedge 2   \big)$,
\begin{equation}\label{e:5.13}
E_{ x }^\vp    \left[ \sup_{0\leq t \leq T}  \big|\wh Y_t^\vp \big|^p    \right] < \infty.
\end{equation}
Meanwhile, from the boundedness of $ \nabla u^0$ and \eqref{e:5.5}, it also follows that
\begin{equation}\label{e:5.14}
E_{ x }^\vp    \left[  \int_0^T  \| \wt Z_r^\vp   \|^2  \, dr
\right] < \infty,
\end{equation}
where we have used the fact that $ \omega \in L^2 (\T ^d)$.

Let $ \tau_n$ be the stopping time
$$
\tau_n :=  \inf\Big\{t\geq 0: | \wh Y_t^\vp | \leq \frac1n \Big\}.
$$
By Ito's formula, we have for all $p\in \big( 1, \frac{ Q(\lambda, d ,q_1, q_2)}{2}   \wedge 2   \big)$,
\begin{eqnarray}
&&\big|  \wh Y_{t \wedge \tau_n }^\vp  \big|^p
+ \frac{p}{2}  \int_{t \wedge \tau_n }^{ T \wedge \tau_n }
\big|  \wh Y_r^\vp  \big|^{p-2}  \,
\big\<    \wh a( X_r^\vp / \vp)  \wt Z_r^\vp,   \wt Z_r^\vp
\big\>  \, dr    \nonumber \\
\label{e:5.15}
&=& \big|  \wh Y_{T \wedge \tau_n }^\vp  \big|^p  + p \int_{t \wedge \tau_n }^{ T \wedge \tau_n }
\big|  \wh Y_r^\vp  \big|^{p-1}  \,
\big(  F^\vp (r, \wt Y_r^\vp  , \wt Z_r^\vp)   - F^\vp (r,0,0)    \big)\, dr     \nonumber \\
&&- p \int_{t \wedge \tau_n }^{ T \wedge \tau_n }
\big|  \wh Y_r^\vp  \big|^{p-1}  \,
\big\< \wt Z_r^\vp , dM_r^\vp   \big\>  +   p \int_{t \wedge \tau_n }^{ T \wedge \tau_n }
\big|  \wh Y_r^\vp  \big|^{p-1}  \,
\big(  \wh h^\vp  (r, \wt Y_r^\vp  )  -    \wh h^\vp  (r, 0)   \big)\, dK_r^\vp.
\end{eqnarray}
In view of $p/2 < 1$ and applying Young's inequality, we obtain
\begin{eqnarray*}
&&E_{ x }^\vp    \left[ \sup_{0\leq t\leq T }   \left|    \int_0^t \big|  \wh Y_r^\vp  \big|^{p-1}  \,
\big\< \wt Z_r^\vp , dM_r^\vp   \big\>
\right|   \right]    \\
&\leq& C E_{ x }^\vp  \left[   \left(    \int_0^T \big|  \wh Y_r^\vp  \big|^{2(p-1)}  \,
\big\<    \wh a( X_r^\vp / \vp)  \wt Z_r^\vp,   \wt Z_r^\vp
\big\>  \, dr
\right)^{\frac12}   \right]     \\
&\leq& C E_{ x }^\vp  \left[  \left(  \sup_{0\leq r \leq T }      \big|  \wh Y_r^\vp  \big|^{p-1}   \right)  \,\cdot
\left( \int_0^T  \| \wt Z_r^\vp \|^2  \, dr \right)^{\frac12}
 \right]    \\
&\leq&   C E_{ x }^\vp  \left[  \sup_{0\leq r \leq T }      \big|  \wh Y_r^\vp  \big|^p     \right]
+ C  E_{ x }^\vp  \left[ \left( \int_0^T  \| \wt Z_r^\vp \|^2  \, dr \right)^{ \frac{p}{2} }
\right]      \\
&\leq&  C E_{ x }^\vp  \left[  \sup_{0\leq r \leq T }      \big|  \wh Y_r^\vp  \big|^p     \right]
+ C  E_{ x }^\vp   \left( \left[ \int_0^T  \| \wt Z_r^\vp \|^2  \, dr  \right]   \right)^{ \frac{p}{2} },
\end{eqnarray*}
which implies the local martingale $\int_0^t   \big|  \wh Y_r^\vp  \big|^{p-1}  \,
\big\< \wt Z_r^\vp , dM_r^\vp   \big\> $ is actually a martingale. Denote
$$ H^\vp(r):= \int_0^r  F^\vp (u,0,0)  \, du
+ \int_0^r   \wh h^\vp  (u, 0)  \, dK_u^{\vp}  ,
$$
and taking expectation in \eqref{e:5.15}, then
\begin{eqnarray*}
&&E_{ x }^\vp    \left[
\big|  \wh Y_{t \wedge \tau_n }^\vp  \big|^p   \right]
+ \frac{p}{2}  E_{ x }^\vp    \left[
  \int_{t \wedge \tau_n }^{ T \wedge \tau_n }
\big|  \wh Y_r^\vp  \big|^{p-2}  \,
\big\<    \wh a( X_r^\vp / \vp)  \wt Z_r^\vp,   \wt Z_r^\vp
\big\>  \, dr   \right]   \\
&\leq&  E_{ x }^\vp    \left[
\big|  \wh Y_{T \wedge \tau_n }^\vp  \big|^p     \right]
+  p C  E_{ x }^\vp    \left[  \int_{t \wedge \tau_n }^{ T \wedge \tau_n }   \big|  \wh Y_r^\vp  \big|^{p-1}  \,
\big( | \wt Y_r^\vp  |   +\| \wt Z_r^\vp  \|  \big) \, dr
\right]    \\
&&+  p C  E_{ x }^\vp    \left[  \int_{t \wedge \tau_n }^{ T \wedge \tau_n }   \big|  \wh Y_r^\vp  \big|^{p-1}  \,
| \wt Y_r^\vp  |  dK_r^\vp
\right]   \\
&\leq&   E_{ x }^\vp    \left[
\big|  \wh Y_{T \wedge \tau_n }^\vp  \big|^p     \right]
+  C  E_{ x }^\vp    \left[  \int_{t \wedge \tau_n }^{ T \wedge \tau_n }   \big|  \wh Y_r^\vp  \big|^{p-1}  \,
\big( | \wh Y_r^\vp  | +  | H^\vp (r)|   +\| \wt Z_r^\vp  \|  \big) \, dr
\right]    \\
&&+   C  E_{ x }^\vp    \left[  \int_{t \wedge \tau_n }^{ T \wedge \tau_n }   \big|  \wh Y_r^\vp  \big|^{p-1}  \,
\big(  | \wh Y_r^\vp  | +  | H^\vp (r)|    \big)  \, dK_r^\vp
\right]  .
\end{eqnarray*}
where we have used the fact that functions $ F^\vp (r,y,z)$  and $ \wh h^\vp  (r, y )$ are uniformly Lipschitz in $ (y,z)$ and $y$, respectively. By Young's inequality, we have
\begin{eqnarray*}
\big|  \wh Y_r^\vp  \big|^{p-1}  \,    | H^\vp (r)|
&\leq& \frac1q  \big|  \wh Y_r^\vp  \big|^{(p-1)q}
+ \frac1p | H^\vp (r)|^p  \\
&=&  \frac1q  \big|  \wh Y_r^\vp  \big|^p
+ \frac1p | H^\vp (r)|^p ,
\end{eqnarray*}
and
\begin{eqnarray*}
\big|  \wh Y_r^\vp  \big|^{p-1}  \,   \| \wt Z_r^\vp  \|
&=& \big|  \wh Y_r^\vp  \big|^{ \frac{p}{2} }  \, \cdot
\left(  \big|  \wh Y_r^\vp  \big|^{ \frac{p}{2} -1 }  \,  \| \wt Z_r^\vp  \|   \right)   \\
&\leq& \frac{1}{ \delta}    \big|  \wh Y_r^\vp  \big|^p
+ \delta \big|  \wh Y_r^\vp  \big|^{p-2}   \,  \| \wt Z_r^\vp  \|^2.
\end{eqnarray*}
Consequently,
\begin{eqnarray*}
&&E_{ x }^\vp    \left[
\big|  \wh Y_{t \wedge \tau_n }^\vp  \big|^p   \right]
+ C (1- \delta)    E_{ x }^\vp    \left[
  \int_{t \wedge \tau_n }^{ T \wedge \tau_n }
\big|  \wh Y_r^\vp  \big|^{p-2}  \,
\| \wt Z_r^\vp  \|^2   \, dr   \right]   \\
&\leq&   E_{ x }^\vp    \left[
\big|  \wh Y_{T \wedge \tau_n }^\vp  \big|^p     \right]
+ CE_{ x }^\vp    \left[
\int_{t \wedge \tau_n }^{ T \wedge \tau_n }  | H^\vp (r)|^p
\, dr  \right]
+ CE_{ x }^\vp    \left[
\int_{t \wedge \tau_n }^{ T \wedge \tau_n }  | H^\vp (r)|^p
\,dK_r^\vp  \right]   \\
&&+   C\Big( 1+ \frac{1}{ \delta}   \Big)   E_{ x }^\vp    \left[
  \int_{t \wedge \tau_n }^{ T \wedge \tau_n }
\big|  \wh Y_r^\vp  \big|^p  \,   dr   \right]
+ C\Big( 1+ \frac{1}{ \delta}   \Big)   E_{ x }^\vp    \left[
  \int_{t \wedge \tau_n }^{ T \wedge \tau_n }
\big|  \wh Y_r^\vp  \big|^p  \,   dK_r^\vp   \right]  .
\end{eqnarray*}
Choosing $\delta$ sufficiently small so that $ 1- \delta >0$ and by the version of Gronwall's lemma in \cite[Lemma 3]{Slomiriski}, we deduce that
\begin{eqnarray*}
E_{ x }^\vp    \left[
\big|  \wh Y_{t \wedge \tau_n }^\vp  \big|^p   \right]
&\leq&  C E_{ x }^\vp    \left[
\big|  \wh Y_{T \wedge \tau_n }^\vp  \big|^p     \right]
+ C\left( E_{ x }^\vp    \left[
\int_0^T   | H^\vp (r)|^p  \,   dr   \right]   +  E_{ x }^\vp   \left[    \int_0^T    | H^\vp (r)|^p  \,   dK_r^\vp   \right]
\right).
\end{eqnarray*}

Let $ \tau_\infty:= \lim_{n\to \infty} \tau_n$, that is $\tau_\infty=\inf\{t\geq 0, \wh Y_t =0  \} $.
For all $ t\in [0,T]$, $   \wh Y_{t \wedge \tau_n }^\vp  $ is dominated by $ \sup_{ 0\leq t \leq T}   |\wh Y_t^\vp |$. By \eqref{e:5.13} and letting $ n $ tends to $\infty$, dominated convergence theorem implies that
\begin{eqnarray*}
E_{ x }^\vp    \left[
\big|  \wh Y_{t \wedge \tau_\infty }^\vp  \big|^p   \right]
&\leq&  C E_{ x }^\vp    \left[
\big|  \wh Y_{T \wedge \tau_\infty }^\vp  \big|^p     \right]
+ C\left( E_{ x }^\vp    \left[
\int_0^T   | H^\vp (r)|^p  \,   dr   \right]   +  E_{ x }^\vp   \left[    \int_0^T    | H^\vp (r)|^p  \,   dK_r^\vp   \right]
\right).
\end{eqnarray*}
In view of
$$
\big| \wh Y_{T \wedge \tau_\infty }^\vp   \big| =
\big| \wh Y_{T  }^\vp   \big| \, 1_{ \{  T \leq \tau_\infty  \} }
\leq \left| g(X_T^\vp)  - g (\wt \omega^\vp (X_T^\vp ) ) + H^\vp  (T)  \right|,
$$
then
\begin{eqnarray}
E_{ x }^\vp    \left[
\big|  \wh Y_{t \wedge \tau_\infty }^\vp  \big|^p   \right]
&\leq&  C E_{ x }^\vp    \big[
\left| g(X_T^\vp)  - g (\wt \omega^\vp (X_T^\vp ) ) + H^\vp (T)  \right|^p     \big]   \nonumber \\
\label{e:5.16}
&&+ C\left( E_{ x }^\vp    \left[
\int_0^T   | H^\vp (r)|^p  \,   dr   \right]   +  E_{ x }^\vp   \left[    \int_0^T    | H^\vp (r)|^p  \,   dK_r^\vp   \right]
\right)    \nonumber\\
&=:& I_{1,\vp}+ C( I_{2,\vp} + I_{3,\vp}).
\end{eqnarray}
Now, we are going to prove $E_{ x }^\vp   [
\big|  \wh Y_{t \wedge \tau_\infty }^\vp  \big|^p   ]$ tends to zero as $\vp \to 0$.
Indeed, by Theorem \ref{T:5.2}, the first term in the right-hand side converges to zero, and for all $ r \in [0,T]$, $E_{ x }^\vp  [| H^\vp (r)|^p ] \stackrel{\vp \to 0}{\longrightarrow}0$. On the other hand, it follows from H\"{o}lder's inequality that
\begin{eqnarray*}
E_{ x }^\vp    [| H^\vp (r)|^p ]  &\leq&
2E_{ x }^\vp  \left[ \left| \int_0^r  F^\vp (u,0,0)  \, du    \right|^p  \right]
+2 E_{ x }^\vp  \left[   \left| \int_0^r   \wh h^\vp  (u, 0)  \, dK_u^{\vp}    \right|^p  \right]    \\
&\leq&  2E_{ x }^\vp  \left[ \left( \int_0^r 1\, dr  \right)^{p-1}   \left( \int_0^r  | F^\vp (u,0,0) |^p \, du    \right)  \right]    \\
&& + 2E_{ x }^\vp  \left[ \left( \int_0^r 1\, dK_u^{\vp}  \right)^{p-1}   \left( \int_0^r  | \wh h^\vp  (u, 0) |^p \, dK_u^{\vp}    \right)  \right]  \\
&\leq&  C T^{p-1}  E_{ x }^\vp  \bigg[  \int_0^T  \big(  1+ \|\nabla\omega^\sharp \|^{2p} \big) \, du      \bigg]   +C E_{ x }^\vp  \bigg[  (K_T^\vp)^{p-1} \,    \Big(\int_0^T \big(  1+ \|\nabla\omega^\sharp \|^p  \big) \, dK_u^{\vp}
 \Big) \bigg]   \\
&\leq& CT^p  \big(  1+ \|\nabla\omega^\sharp \|^{2p} \big)
+ C \big(  1+ \|\nabla\omega^\sharp \|^p \big)   \,
E_{ x }^\vp [ (K_T^{\vp})^{p}] .% < \infty.
\end{eqnarray*}
Hence by dominated convergence, the term $I_{2,\vp}$ in \eqref{e:5.16} converges also to zero.
Since for all $ r \in [0,T]$, $E_{ x }^\vp  [| H^\vp (r)|^p ] \stackrel{\vp \to 0}{\longrightarrow}0$, which implies for any $ \alpha >0$, $ \P_{ x }^\vp  ( | H^\vp (r)|^p > \alpha ) \stackrel{\vp \to 0}{\longrightarrow}0$.
Meanwhile for any fixed $ \delta >0$, we have
\begin{eqnarray*}
&&\P_{ x }^\vp \left[  \left(  \int_0^T    | H^\vp (r)|^p  \,   dK_r^\vp  \right) > \delta   \right]   \\
&=&\P_{ x }^\vp \left[  \left(  \int_0^T    | H^\vp (r)|^p  \,   dK_r^\vp  \right) > \delta  , \, | H^\vp (r)|^p > \alpha \right] \\
&&+ \P_{ x }^\vp \left[  \left(  \int_0^T    | H^\vp (r)|^p  \,   dK_r^\vp  \right) > \delta  , \, | H^\vp (r)|^p \leq  \alpha \right]  \\
&\leq&  \P_{ x }^\vp \left[   | H^\vp (r)|^p > \alpha \right]
+ \P_{ x }^\vp \left[   K_T^\vp   > \delta/\alpha  , \, | H^\vp (r)|^p \leq  \alpha \right].
\end{eqnarray*}
From the definition of $H^\vp $ and the boundedness of $\wh h^\vp  (u, 0) $, it yields for any fixed $ \delta$, we can choose appropriate $ \alpha$ such that the second term in the right hand side equals to $0$. %$ \P_{ x }^\vp \left[   K_T^\vp   > \delta/\alpha  , \, | H^\vp (r)|^p \leq  \alpha \right]=0$.
Thus for any $ \delta>0$, it yields
$$
\P_{ x }^\vp \left[  \left(  \int_0^T    | H^\vp (r)|^p  \,   dK_r^\vp  \right) > \delta   \right]   \stackrel{\vp \to 0}{\longrightarrow}0.   %\qquad \quad \forall   \delta>0 .
$$
Since $ \int_0^T    | H^\vp (r)|^p  \,   dK_r^\vp $ is nonnegative and integrable, then the term $I_{3,\vp}$ in \eqref{e:5.16} converges to zero as $ \vp$ tends to $0$. %$ E_{ x }^\vp  [\int_0^T | H^\vp (r)|^p \,   dK_r^\vp  ] \stackrel{\vp \to 0}{\longrightarrow}0$.
Consequently, we have proved that
$
E_{ x }^\vp    [
\big|  \wh Y_{t \wedge \tau_\infty }^\vp  \big|^p   ]
\stackrel{\vp \to 0}{\longrightarrow}0
$ for all $ t\in [0,T]$.
Taking $t=0$, we have the desired conclusion.
\qed

{\bf Proof of Theorem \ref{T:1.1}.} Let us assume now that $ f$ and $g$ satisfy Assumptions \ref{h:3}-\ref{h:4}.
Let $ \rho_1: \R^d \times \R \times \R^d  \to \R  $ be a $ C_c^\infty$ function with $ \int \rho_1 (x,y,z) dx\, dy\, dz =1$. Define $ f_n := \rho_n * f$ with $ \rho_n (x,y,z)= n^{2d+1} \rho (nx,ny,nz)$, then $f_n $ is infinitely differentiable with bounded derivatives.
Since $ g$ is bounded and continuous, we can also approximate $g$ by a sequence $ g_n$ of functions in $ C_b^\infty (\R^d) $.
%Moreover, for all

Let $ (Y_n^\vp, Z_n^\vp )$ be the solution of the BSDE
\begin{eqnarray*}
Y_n^\vp (s) &=& g_n (X^\vp (T) ) + \int_s^T  f_n( X^\vp (r), Y_n^\vp (r), Z_n^\vp (r)  ) \, dr   + \int_s^T  c( X^\vp (r) / \vp )  Y_n^\vp (r) \, dK_r^\vp     \\
&& - \int_s^T  \< Z_n^\vp (r) , dM^\vp (r)  \>, \quad t\leq s \leq T , \qquad  \P_{t,x}^{\vp} -a.s. ,
\end{eqnarray*}
satisfying
$
E_{ t,x }^\vp  \left[  \sup_{t\leq s \leq T}   | Y_n^\vp (s)  |^2   + \int_t^T  \| Z_n^\vp (s)  \|^2 ds
\right] < \infty.
$
By Ito's formula, we have
\begin{eqnarray*}
&& e^{ 2\int_0^s c_1^+ ( X_r^\vp ) dr }   | Y_n^\vp (s)- Y^\vp (s)  |^2  + \int_s^T  e^{ 2\int_0^r c_1^+ (X_u^\vp) du }  \big\<  a( X(r))  ( Z_n^\vp (r)- Z^\vp (r)  ), Z_n^\vp (r)- Z^\vp (r)  \big\>   dr   \\
&=&  e^{ 2\int_0^T c_1^+  (X_r^\vp) dr }  | g_n (X^\vp (T) )- g (X^\vp (T) )  |^2  - 2\int_s^T  c_1^+ (X_r^\vp)  e^{ 2\int_0^r c_1^+ (X_u^\vp) du }\,
 | Y_n^\vp (r)- Y^\vp (r)  |^2 \, dr
\\
&&- 2 \int_s^T    e^{ 2\int_0^r c_1^+ (X_u^\vp) du }
\big( Y_n^\vp (r)- Y^\vp (r)  \big)\,
\big\<  Z_n^\vp (r)- Z^\vp (r)  , dM^\vp (r)   \big\>
\\
&&+ 2 \int_s^T    e^{ 2\int_0^r c_1^+ (X_u^\vp) du }
\big( Y_n^\vp (r)- Y^\vp (r) \big)\,   \big[f_n( X^\vp (r), Y_n^\vp (r), Z_n^\vp (r)  ) -  f( X^\vp (r), Y^\vp (r), Z^\vp (r)  )  \big]  \, dr   \\
&&+ 2 \int_s^T  e^{ 2\int_0^r c_1^+ (X_u^\vp) du }  c( X^\vp (r) / \vp )   | Y_n^\vp (r)- Y^\vp (r) |^2   \, dK_r^\vp  .
\end{eqnarray*}
Since function $ c$ is nonpositive,
$ - ( c_1^+ (X_r^\vp) - c_1 (X_r^\vp) ) \leq 0$
and
\begin{eqnarray*}
&& f_n( X^\vp (r), Y_n^\vp (r), Z_n^\vp (r)  ) -  f( X^\vp (r), Y^\vp (r), Z^\vp (r)  )   \\
&=& [f_n( X^\vp (r), Y_n^\vp (r), Z_n^\vp (r)  )  - f_n( X^\vp (r), Y^\vp (r), Z_n^\vp (r)  )  ]  \\
&&+[ f_n( X^\vp (r), Y^\vp (r), Z_n^\vp (r)  )  -  f_n( X^\vp (r), Y^\vp (r), Z^\vp (r)  )  ]  \\
&&+[  f_n( X^\vp (r), Y^\vp (r), Z^\vp (r)  )   -  f( X^\vp (r), Y^\vp (r), Z^\vp (r)],
\end{eqnarray*}
we can obtain by usual computations as \eqref{e:4.5} that
\begin{eqnarray}\label{e:5.19}
&& e^{ 2\int_0^s c_1^+ ( X_r^\vp ) dr }   | Y_n^\vp (s)- Y^\vp (s)  |^2  + \int_s^T  e^{ 2\int_0^r c_1^+ (X_u^\vp) du }  \big\<  a( X(r))  ( Z_n^\vp (r)- Z^\vp (r)  ), Z_n^\vp (r)- Z^\vp (r)  \big\>   dr    \nonumber \\
&\leq& e^{ 2\int_0^T c_1^+  (X_r^\vp) dr }  | g_n (X^\vp (T) )- g (X^\vp (T) )  |^2    \nonumber \\
&&- 2 \int_s^T    e^{ 2\int_0^r c_1^+ (X_u^\vp) du }     \big( Y_n^\vp (r)- Y^\vp (r)  \big)\,
\big\<  Z_n^\vp (r)- Z^\vp (r)  , dM^\vp (r)   \big\>
 \nonumber \\
&&+ 2c_2 \int_s^T    e^{ 2\int_0^r c_1^+ (X_u^\vp) du }     \big| Y_n^\vp (r)- Y^\vp (r)  \big|\,
\big| Z_n^\vp (r)- Z^\vp (r) ) \big|\, dr
 \nonumber \\
&&+ 2  \int_s^T    e^{ 2\int_0^r c_1^+ (X_u^\vp) du }     \big| Y_n^\vp (r)- Y^\vp (r)  \big|\,
\big| (f_n - f) (X^\vp (r), Y^\vp (r), Z^\vp (r) ) ) \big|\, dr
 \nonumber \\
&\leq& e^{ 2\int_0^T c_1^+  (X_r^\vp) dr }  | g_n (X^\vp (T) )- g (X^\vp (T) )  |^2  + C_\lambda \int_s^T    e^{ 2\int_0^r c_1^+ (X_u^\vp) du }     \big| Y_n^\vp (r)- Y^\vp (r)  \big|^2  \, dr
 \nonumber \\
&&+\frac12  \int_s^T  e^{ 2\int_0^r c_1^+ (X_u^\vp) du }  \big\<  a( X(r))  ( Z_n^\vp (r)- Z^\vp (r)  ), Z_n^\vp (r)- Z^\vp (r)  \big\>   dr
 \nonumber \\
&&+ \int_s^T    e^{ 2\int_0^r c_1^+ (X_u^\vp) du }
\big| (f_n - f) (X^\vp (r), Y^\vp (r), Z^\vp (r) ) ) \big|^2\, dr    \nonumber \\
&&- 2 \int_s^T    e^{ 2\int_0^r c_1^+ (X_u^\vp) du }     \big( Y_n^\vp (r)- Y^\vp (r)  \big)\,
\big\<  Z_n^\vp (r)- Z^\vp (r)  , dM^\vp (r)   \big\>.
\end{eqnarray}
Taking expectation and by Gronwall's inequality, it yields
\begin{eqnarray}\label{e:5.20}
&& E_{ x }^\vp \big[e^{ 2\int_0^s c_1^+ ( X_r^\vp ) dr }   | Y_n^\vp (s)- Y^\vp (s)  |^2   \big]   \nonumber \\
&&+ E_{ x }^\vp   \Big[\int_s^T  e^{ 2\int_0^r c_1^+ (X_u^\vp) du }  \big\<  a( X(r))  ( Z_n^\vp (r)- Z^\vp (r)  ), Z_n^\vp (r)- Z^\vp (r)  \big\>   dr  \Big]    \nonumber \\
&\leq& C E_{ x }^\vp    \bigg[e^{ 2\int_0^T  c_1^+  (X_r^\vp) dr }  | g_n (X_T^\vp )- g (X_T^\vp )  |^2    \nonumber \\
&&  +  \int_s^T  e^{ 2\int_0^r   c_1^+  (X_u^\vp) du }\, | f_n( X_r^\vp , Y_r^\vp , Z_r^\vp   ) - f( X_r^\vp , Y_r^\vp , Z_r^\vp   ) |^2  \, dr    \bigg] .
\end{eqnarray}
By both Burkh\"{o}lder's and Young's inequalities, it implies
\begin{eqnarray}\label{e:5.21}
&&2 E_{ x }^\vp  \bigg[ \sup_{s\in [0,T]}  \bigg|   \int_0^s    e^{ 2\int_0^r c_1^+ (X_u^\vp) du }     \big( Y_n^\vp (r)- Y^\vp (r)  \big)\,  \big\<  Z_n^\vp (r)- Z^\vp (r)  , dM^\vp (r)   \big\>
\bigg|     \bigg]   \nonumber \\
&\leq& C E_{ x }^\vp \bigg[ \bigg( \int_0^T  e^{ 4\int_0^r c_1^+ (X_u^\vp) du }  \big|Y_n^\vp (r)- Y^\vp (r)  \big|^2 \,
 \big\<  a( X(r))  ( Z_n^\vp (r)- Z^\vp (r)  ), Z_n^\vp (r)- Z^\vp (r)  \big\>   dr
\bigg)^{\frac12}   \bigg]    \nonumber \\
&\leq& \frac12   E_{ x }^\vp    \Big[ \sup_{s\in [0,T]}  e^{ 2\int_0^s  c_1^+ (X_r^\vp) dr } | Y_n^\vp (s)- Y^\vp (s)  |^2  \Big]      \nonumber \\
&&\quad + C  E_{ x }^\vp   \Big[\int_s^T  e^{ 2\int_0^r c_1^+ (X_u^\vp) du }  \big\<  a( X(r))  ( Z_n^\vp (r)- Z^\vp (r)  ), Z_n^\vp (r)- Z^\vp (r)  \big\>   dr  \Big].
\end{eqnarray}
Consequently, in view of \eqref{e:5.19}, \eqref{e:5.20} and \eqref{e:5.21}, we conclude that
\begin{eqnarray}
&&E_{ x }^\vp    \Big[ \sup_{s\in [0,T]} | Y_n^\vp (s)- Y^\vp (s)  |^2  \Big]   \nonumber \\
\label{e:5.17}
&\leq& E_{ x }^\vp    \Big[ \sup_{s\in [0,T]}  e^{ 2\int_0^s  c_1^+ (X_r^\vp) dr } | Y_n^\vp (s)- Y^\vp (s)  |^2  \Big]
\nonumber \\
&\leq& C E_{ x }^\vp    \bigg[e^{ 2\int_0^T  c_1^+  (X_r^\vp) dr }  | g_n (X_T^\vp )- g (X_T^\vp )  |^2   \nonumber \\
&&\quad  +  \int_0^T  e^{ 2\int_0^r   c_1^+  (X_u^\vp) du }\, | f_n( X_r^\vp , Y_r^\vp , Z_r^\vp   ) - f( X_r^\vp , Y_r^\vp , Z_r^\vp   ) |^2  \, dr
 \bigg]   \nonumber \\
&\leq&  C \Big( \| g_n - g  \|_{L^\infty }^2    E_{ x }^\vp \Big[ e^{ 2\int_0^T  c_1^+  (X_r^\vp) dr } \Big]    + \omega_n (f)^2  \,
E_{ x }^\vp \Big[  \int_0^T  e^{ 2\int_0^r  c_1^+ (X_u^\vp) du }\, dr \Big]     \Big)     \nonumber \\
&\leq&  C \Big( \| g_n - g  \|_{L^\infty }^2       + \omega_n (f)^2      \Big)   ,
\end{eqnarray}
where the fact for all $ (x,y,z)\in\R^d \times \R \times \R^{d}$,
\begin{eqnarray*}
| f_n (x,y,z)- f(x,y,z)  | \leq \sup_{ \| (x,y,z)- (x',y',z') \| \leq   \frac{1}{n}  }   |f(x,y,z)- f (x',y',z') |  =: \omega_n (f),
\end{eqnarray*}
has been used.

In the same way, define $ \bar{f_n} (x,y,z):= \int_{\T^d} f_n \big(x,y,z \, \nabla\wt \omega ( \eta)   \big)\, m(\eta) d\eta$. Then $ \bar{f_n}$ satisfies Assumption \ref{h:3} with constant independent of $n$. Also for all $ (x,y,z)$, we have
$ | \bar{f}_n (x,y,z)-  \bar{f} (x,y,z)  |\leq \omega_n (f)$.
Let $E_x^0  $ denotes the expectation under the law of a reflected Brownian motion $ X^0$ with covariance matrix $ \bar {a}$, $ M^0 $ and $ K^0$ be the martingale part and local time of $ X^0$.
Then if $( \bar{Y}_n, \bar{Z}_n)$ is the solution of the BSDE
\begin{eqnarray*}
\bar{Y}_n (t) &=& g_n (X^0 (T) ) + \int_t^T  \bar{f_n} \big(X_r^0 ,\bar{Y}_n (r), \bar{Z}_n (r) \big)  \, dr   \\
&&-\int_t^T   \big\<\bar{Z}_n (r), dM^0 (r)   \big\>    + \bar{C} \int_t^T   \bar{Y}_n (r)  \, dK_r^0,
\quad t\in [0, T] , \qquad  \P_x^0 -a.s. ,
\end{eqnarray*}
and if $( \bar{Y}, \bar{Z})$ is the solution of the BSDE
\begin{eqnarray*}
\bar{Y} (t) &=& g (X^0 (T) ) + \int_t^T  \bar{f} \big(X_r^0 ,\bar{Y} (r), \bar{Z} (r) \big)  \, dr   \\
&&-\int_t^T   \big\<\bar{Z} (r), dM^0 (r)   \big\>    + \bar{C} \int_t^T   \bar{Y} (r)  \, dK_r^0,
\quad t\in [0, T] , \qquad  \P_x^0 -a.s. ,
\end{eqnarray*}
we similarly have
\begin{eqnarray}
&&E_{ x }^0    \Big[ \sup_{s\in [0,T]} | \bar{Y}_n  (s)- \bar{Y} (s)  |^2  \Big]  \nonumber \\
\label{e:5.18}
&\leq& C \Big( \| g_n - g  \|_{L^\infty }^2   E_{ x }^0 \Big[ e^{ 2\int_0^T  c_1^+  (X_r^0) dr } \Big]    + \omega_n (f)^2  \,
E_{ x }^0 \Big[  \int_0^T  e^{ 2\int_0^r  c_1^+ (X_u^0) du }\, dr \Big]  \Big)     \nonumber \\
&\leq&  C \Big( \| g_n - g  \|_{L^\infty }^2       + \omega_n (f)^2      \Big)   .
\end{eqnarray}
From \eqref{e:5.17}, \eqref{e:5.18} and Lemma \ref{L:5.4}, we can conclude that $ \lim_{\vp \to 0} | Y^\vp (0)- Y(0)  |^p =0$ holds for any $ p\in \big( 1, \frac{ Q(\lambda, d ,q_1, q_2)}{2}   \wedge 2   \big)$.

%\section{references references}

\end{document}